\documentclass[twoside]{article}

\usepackage{preamble}

\usepackage[accepted]{aistats2020}
\begin{document}

%

%

\twocolumn[

\aistatstitle{Nonparametric Estimation in the Dynamic Bradley-Terry Model}

\aistatsauthor{Heejong Bong* \And Wanshan Li* \And  Shamindra Shrotriya* \And Alessandro Rinaldo}


\aistatsaddress{ Department of Statistics and Data Science \\
Carnegie Mellon University } 
]

\begin{abstract}
  We propose a time-varying generalization of the Bradley-Terry model that allows for nonparametric modeling of dynamic global rankings of distinct teams. We develop a novel estimator that relies on kernel smoothing to pre-process the pairwise comparisons over time and is applicable in sparse settings where the Bradley-Terry may not be fit. We obtain  necessary and sufficient conditions for the existence and uniqueness of our estimator. We also derive time-varying oracle bounds for both the estimation error and the excess risk in the model-agnostic setting where the Bradley-Terry model is not necessarily the true data generating process.  We thoroughly test the practical effectiveness of our model using both simulated and real world data and suggest an efficient data-driven approach for bandwidth tuning.
\end{abstract}

\section{Introduction and Prior Work}
\subsection{Pairwise Comparison Data and the Bradley-Terry Model}
Pairwise comparison data are very common in daily life, especially in cases where the goal is to rank several objects. Rather than directly ranking all objects simultaneously, it is usually much easier and more efficient to first obtain results of pairwise comparisons and then use them to derive a \textit{global} ranking across all individuals in
a principled manner. Since the required global rankings are not directly observable, developing a statistical framework for the estimating rankings is a challenging unsupervised learning problem.
One such statistical model for deriving global rankings using pairwise comparisons was presented in the classic paper \citep{bradley1952rank}, and thereafter commonly referred to as the Bradley-Terry model in the literature. A similar model was also studied by Zermelo \citep{Zermelo1929}. The Bradley-Terry model is one of the most popular models to analyze paired comparison data due to its interpretability and computational efficiency in parameter estimation. The Bradley-Terry model along with its variants has been studied and applied in various ranking applications across many domains. This includes the ranking of sports teams \citep{MaV2012, CMV2012, FaT1994}, scientific journals \citep{St1994, Va2016}, and the quality of several brands \citep{Ag2002, RaJ2007}, to name a few. 
\par
In order to introduce the  Bradley-Terry model, suppose that we have $N$ distinct teams, each with a positive
score $s_i$, $i \in [N] \defined \theset{1, \ldots, N}$, quantifying its propensity to be picked or win over other items. The model postulates that the comparisons between different pairs are independent and the results of the comparisons between a given pair, say team $i$ and team $j$, are independent and identically distributed Bernoulli random variables, with {\it winning probability} defined as

\begin{equation}\label{eqn:bradley_terry_prob_succ}
    \prob{i\ \text{ beats } j} = \frac{s_i}{s_i + s_j}, \: \forall \; i,j \in [N]
\end{equation}
A common way to parametrize the model is to set, for each $i$,  $s_i = \exp(\beta_i)$, where $(\beta_1,\ldots,\beta_N)$ are real parameters such that $\sum_{i \in [N]} \beta_i = 0$ (this latter condition is  to ensure identifiability). In this case, equation \eqref{eqn:bradley_terry_prob_succ} is usually expressed as
$\logit(\prob{i\ \text{ beats } j}) = \beta_i - \beta_j$, where, for $x \in (0,1)$, $\logit(x) \defined \log\frac{x}{1-x}$.

\subsection{The time-varying (dynamic) Bradley-Terry Model}
In many applications it is very common to observe paired comparison data spanning over multiple (discrete) time periods. A natural question of interest is then to understand how the \textit{global} rankings \textit{vary} over time i.e. \textit{dynamically}. For example, in sports analytics the performance of teams often changes across match rounds and thus explicitly incorporating the time-varying dependence into the model is crucial. In particular the paper \citep{FaT1994} considers a state-space generalization of the Bradley-Terry model to analyze sports tournaments data. In a similar manner, bayesian frameworks for the dynamic Bradley-Terry model are studied further in \citep{Gli1993, glickman1998statespace_nfl, lopez2018often}. Such dynamic ranking analysis is becoming increasingly important because of the rapid growth of openly available time-dependent paired comparison data.
\par

Our main focus in this paper is to tackle the problem of generalizing the Bradley-Terry model to the time-varying setting with statistical guarantees. Our approach estimates the \textit{changes} in the Bradley-Terry model parameters over time nonparametrically. This enables the derivation of time-varying global dynamic rankings with minimal parametric assumptions. Unlike previous time-varying Bradley-Terry estimation approaches, we seek to establish guarantees in the model-agnostic setting where the Bradley-Terry model is not the true data generating process. This is in contrast to more assumption-heavy parametric frequentist dynamic Bradley-Terry models including \citep{CMV2012}. Our method is also computationally efficient compared to some state-space methods including but not limited to \citep{FaT1994,glickman1998statespace_nfl,MKG2019}.

\section{Time-varying Bradley-Terry Model}\label{sec:our_proposed_bttv}

\subsection{Model Setup}

In our time-varying generalization of the original Bradley-Terry model we assume that $N$ distinct teams play against each other at possibly different times, over a given time interval which, without loss of generality, is taken to be $[0, 1]$. The result of a game between team $i$ and team $j$ at time $t$ is determined by the timestamped winning probability $p_{ij}(t) \defined \prob{i \text{ defeats } j \text{ at time } t}$ which we assume arises from a distinct Bradley-Terry model, one for each time point. In detail, for each $i$, $j$, and $t$
\begin{equation}\label{eqn:bradley_terry_prob_succ_timevarying}
    \logit(p_{ij}(t)) = \beta_i(t) - \beta_j(t), \: \forall \; i,j \in [N], t \in [0, 1]
\end{equation}
where $\bbeta(t) = \text{vec}(\beta_1(t), \beta_2(t), \dots, \beta_N(t)) \in \reals^N$ is an unknown vector such that $\sum_i \beta_i(t) = 0$. 


We observe the outcomes of $M$ timestamped pairwise matches among the $N$ teams $\theset{(i_m, j_m, t_m) : m \in [M]}$. Here $(i_m, j_m, t_m)$ indicates that team $i_m$ and team $j_m$ played at time $t_m$, where $t_1 \le t_2 \le \ldots \le t_M$. The result of the $m$-th match can be expressed in a $N \times N$ data matrix $X^{(m)}$ as follows:
\begin{equation}
    X_{ij}^{(m)} = \begin{cases}
    \begin{split}
        & \mathbf{1}(i \text{ defeats } j \text{ at } t_m) \\
        & \sim \text{Bernoulli}(p_{ij}(t_m))
    \end{split}
    & \text{for } i = i_m \text{ and } j = j_m \\
    1 - X_{ji}^{(m)} & \text{for } i = j_m \text{ and } j = i_m \\
    0 & \text{elsewhere}
    \end{cases}
\end{equation}

Our goal is to estimate
the underlying parameters $\bbeta(t)$ where $t \in [0, 1]$, and then derive the corresponding global ranking of teams. In order to make the estimation problem tractable we assume that the  parameters $\bbeta(t)$ vary smoothly as a function of time $t \in [0, 1]$. 
It is worth noting  that the naive strategy of estimating the model parameters separately at each observed discrete time point on the original data will suffer from two major drawbacks: (i) it will in general  not guarantee smoothly varying estimates and, perhaps more importantly, (ii) computing the maximum likelihood estimator (MLE) of the  parameters in the static Bradley-Terry model may be infeasible due to sparsity in 
the  data (e.g., at each time point we may observe only one match), as we discuss below in \ref{sec:existence_uniqueness_soln}. To overcome these issues we propose a nonparametric methodology which involves kernel-smoothing
the observed data over time.

\subsection{Estimation}\label{sec:estimator}
Our approach in estimating time-varying global rankings is described in the following three step procedure:

\begin{enumerate}
    \item \textit{Data pre-processing:} Kernel smooth the pairwise comparison data across all time periods. This is used to obtain the smoothed pairwise data at each time $t$:
    \begin{equation}\label{eqn:smoothed_data}
        \tilde{X}(t) = \sum_{m=1}^M W_h(t_m, t) X^{(m)},
    \end{equation}
    where $W_h$ is an appropriate kernel function  with bandwidth $h$, which controls the extent of data smoothing. The higher the value of $h$ is, the smoother $\tilde{X}_{ij}(t)$ is over time. 
    \item \textit{Model fitting:} Fit the regular Bradley-Terry
    model on the smoothed data $\tilde{X}_{ij}(t)$. The model estimates the performance of each team at time $t$ using the estimated score vector 
    \begin{equation}\label{eqn:opt_problem}
        \widehat{\bbeta}(t) = {\arg\min}_{\bbeta: \sum_i \beta_{i} (t) = 0} \widehat{\mathcal{R}}(\bbeta;t)
    \end{equation}
     minimizing the negative log-likelihood risk
    \begin{equation}\label{eqn:emp_risk}
        \widehat{\mathcal{R}}(\bbeta;t) = \sum_{i,j: i \neq j} \frac{\tilde{X}_{ij}(t)\log(1 + \exp(\beta_{j}(t) - \beta_{i}(t)))}{\sum_{i,j: i \neq j} \tilde{X}_{ij}(t)} 
    \end{equation}

    \item \textit{Derive global rankings}: Rank each team at time $t$ by its score from $\hat{\bbeta}(t)$.
\end{enumerate}

We observe that if $t = t_1 = t_2 = \dots = t_M$ then step 1 reduces to the original (static) Bradley-Terry model. In this case,
\begin{equation}
\begin{split}
    \tilde{X}_{ij}(t) &= W_h(t,t) \sum_{m = 1}^{M} \mathbf{1}(i_m = i, j_m = j) \\
                      &= W_h(t,t) X_{ij}
\end{split}
\end{equation}
where $X_{ij} = \#\{i \text{ defeated } j\}$. Thus, fitting the model on $\tilde{X}(t)$ in step 2 is equivalent to the original method on data $X$. In this sense our proposed ranking approach represents a time-varying generalization of the original Bradley-Terry model.


This \textit{data pre-processing} is a critical step in our method and is similar to the approach adopted in \citep{zhou10_time_varyin_undir_graph} where it was used in the context of estimating smoothed time varying undirected graphs. This approach has two main advantages. First, applying kernel smoothing on the input pairwise comparison data enables borrowing of information across timepoints. In sparse settings, this reduces the data requirement at each time point to meet necessary and sufficient conditions required for the Bradley-Terry model to have a unique solution as detailed in Section \ref{sec:existence_uniqueness_soln}. Second, kernel smoothing is computationally efficient in a single dimension and is readily available in open source scientific libraries.

\section{Our Contributions}\label{sec:our_contrib}

Our main contributions in this paper are summarized as follows:

\begin{enumerate}
    \item We obtain necessary and sufficient conditions for the existence and uniqueness for our time-varying estimator $\widehat{\bbeta}(t)$ for each $t \in [0, 1]$ simulataneously. See~\cref{sec:existence_uniqueness_soln}.
    \item We extend the results of \cite{simons1999asymptotics} and obtain statistical guarantees for our proposed method in the form of convergence results of the estimated model parameters \textit{uniformly} over all times. We express such guarantees in the form of oracle inequalities in the model-agnostic setting where the Bradley-Terry model is not necessarily the true data generating process. See~\cref{sec:statistical_property}.
    \item We apply our estimator with an data-driven tuned (by LOOCV) hyperparameter successfully to simulations and to real life applications including a comparison to 5 seasons of NFL ELO ratings. See~\cref{sec:Experiments} and \cref{sec:application_nfl}. 
\end{enumerate}

\section{Existence and uniqueness of solution}\label{sec:existence_uniqueness_soln}
The existence and uniqueness of solutions for model \eqref{eqn:opt_problem} is not guaranteed in general. This is an innate property of the original Bradley-Terry model \citep{bradley1952rank}. As pointed out in \citep{ford1957solution} existence of the MLE for the Bradley-Terry model parameters demands a sufficient amount of pairwise comparison data so that there is enough information of relative performance between any pair of two teams for parameter estimation purposes. For example, if there is a team which has never been compared to the others, there is no information which the model can exploit to assign a score for the team. As such its derived rank could be arbitrary. In addition if there are several teams which have never outperformed the others then the Bradley-Terry model would assign negative infinity for the performance of these teams. It would not be possible to compare amongst them for global ranking purposes. In all such cases, the model parameters are not estimable.

 \cite{ford1957solution} derived the necessary and sufficient condition for the existence and uniqueness of the MLE in the original Bradley-Terry model. Below we show how this condition can also be adapted to guarantee the existence and uniqueness of the solution in our time-varying Bradley-Terry model. The  condition can be stated for each time $t$ in terms of the corresponding kernel-smoothed data $\tilde{X}(t)$ as follows:

\begin{ncond}\label{cond:nec_suff_bt_1} 
In every possible partition of the teams into two nonempty subsets, some team $i$ in the first set and some team $j$ in the second set satisfy $\tilde{X}_{ij}(t) > 0$. Or equivalently, for each ordered pair $(i, j)$, there exists a sequence of indices $i_0=i, i_1, \dots, i_n=j$ such that $\tilde{X}_{i_{k-1} i_k}(t) > 0$ for $k=1,\dots,n$.
\end{ncond}

\begin{nrmk}
If we regard $[|\tilde{X}(t)|_{ij}]$ as the adjacency matrix of a (weighted) directed graph, then \cref{cond:nec_suff_bt_1} is equivalent to the \textit{strong connectivity} of the graph.
\end{nrmk}

Under condition \ref{cond:nec_suff_bt_1} we obtain the following existence and uniqueness theorem on the solution set of the time-varying Bradley-Terry model.

\begin{nthm}\label{thm:exist_unq_thm}
    If the smoothed data $\tilde{X}(t)$ satisfies Condition~\ref{cond:nec_suff_bt_1}, then the solution of \eqref{eqn:opt_problem} uniquely exists at time $t$.
\end{nthm}

Hence, in the proposed time-varying Bradley-Terry model we do not require the strong conditions of \citep{ford1957solution} to be met at each time point, but simply require the \textit{aggregated conditions} in Theorem \ref{thm:exist_unq_thm} to hold. This is a significant weakening of the data requirement. For example, even if one team did not play any game in a match round -- a situation that would preclude the MLE in the standard Bradley-Terry model -- it is still possible to assign a rank to this team in their missing round,  as long a game is recorded in another round (with at least one win and one loss). In this sense, the kernel-smoothing of the data in our time-varying Bradley-Terry model reduces the required sample complexity for a unique solution.

\section{Statistical Properties of the Time-varying Bradley-Terry Model}\label{sec:statistical_property}

\subsection{Preliminaries}
\label{sec:prelim}
Existing results \citep{simons1999asymptotics, negahban2016rank} demonstrate the consistency of the estimated \textit{static} Bradley-Terry scores provided that the data were generated from the Bradley-Terry model. However this assumption may be too restrictive for data generation processes in real world applications. In the rest of this section, we will consider  \textit{model-agnostic} time-varying settings where the Bradley-Terry model is not necessarily the true pairwise data generating model. In order to investigate the statistical properties of the proposed estimator, we impose the following relatively mild assumptions.

\begin{nassmp}\label{assmp:regularity}
    Each pair of teams $(i,j)$ play $T^{(i,j)}$ times at time points $\{ t_k^{(i,j)}, k = 1, 2, \dots, T^{(i,j)} \}$, where each  $T^{(i,j)}>0$ satisfy the following conditions, for fixed constants $T>0$ and $0< D_m \leq 1 \leq D_M$:
    \begin{enumerate}
        \item $T^{(i,j)} > T$ ;
        \item for every interval $(a,b) \subset [0, 1]$,
        \begin{equation}\label{eqn:bounded_density}
        \begin{split}
             \lfloor D_m (b- a) T^{(i, j)} \rfloor \leq& \abs{\theset{k: t_k^{(i,j)} \in (a, b)}} \\
             \leq& \lceil D_M (b-a) T^{(i,j)} \rceil. 
        \end{split}
        \end{equation}
    \end{enumerate}
   \end{nassmp}
   
     We remark that the second condition further implies that
    \begin{equation}\label{eqn:bounded_distance}
    \begin{split}
        & t_1^{(i,j)} \leq \frac{1}{D_m T^{(i,j)}}, \enspace t_{T^{(i,j)}}^{(i,j)} \geq 1 - \frac{1}{D_m T^{(i,j)}},\\
        & \frac{1}{D_M T^{(i,j)}} \leq t_{k+1}^{(i,j)} - t_{k}^{(i,j)} \leq \frac{1}{D_m T^{(i,j)}}
    \end{split}
    \end{equation} 
    for $k = 1, 2, \dots, T^{(i,j)}-1$.

\cref{assmp:regularity} allows for different team pairs to play against each other a different number of times and for the game times to be spaced irregularly, though in a controlled manner. To enable statistical analyses on time-varying quantities, we further require that the winning probabilities satisfy a minimal degree of smoothness and that their rate of decay is controlled.

\begin{nassmp}\label{assmp:uniform_smooth_bound}
    For any $i,j$, the function $t \in [0,1] \mapsto p_{ij}(t)$ is Lipschitz with universal constant $L_p$ and uniformly bounded below $p_{\min} > 0$ which is dependent to $N$ and $T$.
\end{nassmp}
The quantity $p_{\min}$ need not be bounded away from $0$ as function of $T$ and $N$. However, in order to guarantee estimability of the model parameters in time-varying settings,  we will need to control the rate at which it is allowed to vanish. See \cref{thm:asymp_as_cond} below.

Finally, we assume that the kernel used to smooth over time satisfy the following regularity conditions, which are quite standard in the nonparametric literature.

\begin{nassmp}\label{assmp:kernel_function}
The kernel function 
     $W:(-\infty, \infty) \rightarrow (0, \infty)$ is a symmetric function such that
    \begin{equation}
    \begin{split}
        \int_{-\infty}^{\infty} W(x) dx = 1, &\enspace \int_{-\infty}^{\infty} |x| W(x) dx < \infty \\
        \mathcal{V}(W) < \infty, &\enspace \mathcal{V}(|\cdot|W) < \infty \\
    \end{split}
    \end{equation}
    where $\mathcal{V}(f(x))$ is the total variation of a function $f$.
    For each $s,t \in [0,1]$ we further write 
    \begin{equation}
        W_h(s,t) = \frac{1}{h} W\left(\frac{s-t}{h}\right).
    \end{equation}
\end{nassmp}

It is easy to see that these conditions imply 
\begin{equation}
    \|W\|_\infty = \sup_x W(x) < \infty.
\end{equation}
Thus, without loss of generality, we assume that $\|W\|_\infty \leq 1$; the general case can be handled by modifying the constants accordingly. The use of kernels satisfying the above assumptions is standard in nonparametric problems involving H\"older continuous functions of order $1$, such as the winning probabilities function of Assumption~\ref{assmp:uniform_smooth_bound}.

\subsection{Existence and uniqueness of solutions}

 \cite{simons1999asymptotics} showed that the necessary and sufficient condition for the existence and uniqueness of the MLE in the original Bradley-Terry model is satisfied asymptotically almost surely under minimal assumptions. Below, we show that this type of result can be extended to our more general time-varying settings.
\begin{nthm}\label{thm:asymp_as_cond}
    \begin{equation}
    \begin{split}
        & \Pr(\text{Condition~\ref{cond:nec_suff_bt_1} is satisfied at every } t)  \\
        & \geq 1 - 4N\exp\left( - \frac{NT}{2}  p_\text{min} \right).
    \end{split}
    \end{equation}
\end{nthm}

\begin{nrmk}
As we remarked above, $p_{\min}$ needs not be bounded away from zero, but is allowed to vanish slowly enough in relation to $N$ and $T$ so that condition \eqref{thm:asymp_as_cond} is fulfilled as long as $\frac{1}{p_\text{min}} = o\left(\frac{NT}{2\log N}\right)$.
\end{nrmk}


\subsection{Oracle Properties}

In our general agnostic time-varying setting the Bradley-Terry model is not assumed to be the true data generating process. It follows that, for each $t$, there is no true parameter to which to compare the estimator defined in \eqref{eqn:opt_problem}. Instead, we may compare it to the \textit{projection} parameter $\bbeta^*(t) \in \mathbb{R}^N$, which is the best approximation to the winning probabilities at time $t$ using the dynamic Bradley Terry model; see \eqref{eqn:bradley_terry_prob_succ_timevarying}. In detail, the oracle parameter is defined as 
\begin{equation}
    \bbeta^*(t) = {\arg\min}_{\bbeta: \sum_i \beta_{i}(t) = 0} \mathcal{R}(\bbeta; t)
\end{equation}
where
\begin{equation}\label{eqn:oracle_risk}
    \mathcal{R}(\bbeta; t) = \frac{1}{\binom{N}{2}} \sum_{i,j: i \neq j} p_{ij}(t) \log(1 + \exp(\beta_{j}(t) - \beta_{i}(t)))
\end{equation}
We note that, when the winning probabilities obey a Bradley-Terry model,  the projection parameter corresponds to the 
true model parameters.

Next, for each fixed time $t \in (0,1)$ and $h > 0$, we introduce two quantities, namely $M(t)$ and $\delta_h(t)$, that can be thought of as conditions numbers of sort, affecting directly both the estimation and prediction accuracy of the proposed estimator. In detail, we set
\[
    M(t) = \max_{i,j: i \neq j}\exp(\beta_i^*(t) - \beta_j^*(t)) \]
    and
    \[
    \delta_h(t) = \max_i \sum_{j: j \neq i} \left| \frac{\tilde{T}_{ij}(t)}{\tilde{T}_i(t)} - \frac{1}{N-1} \right|.
\]
where $\tilde{T}_{ij}(t) = \tilde{X}_{ij}(t) + \tilde{X}_{ji}(t)$ and $\tilde{T}_i(t) = \sum_{j: j \neq i} \tilde{T}_{ij}(t)$. The ratio $M(t)$ quantifies the maximal discrepancy in winning scores among all possible pairs at time $t$, and, as shown in \cite{simons1999asymptotics}, determines the consistency rate of the MLE in the traditional Bradely-Terry model (see also the comments following \cref{remark:compare.to.Simons.Yao} below). The quantity $\delta_h(t)$ is instead a measure of regularity in how evenly the teams play against each other. In particular $\delta_h(t) = 0$, for all $t$ and $h$ when there is a constant number of matches among each pair of teams, for each time. Since we allow for the possibility of different number of matches between teams and across time, it becomes necessary to quantity such degree of design irregularity.

In order to verify the quality of the proposed estimator $\widehat{\bbeta}(t)$, we will consider high-probability oracle bound on estimation error $\norm{\widehat{\bbeta}(t) - \bbeta^*(t)}_{\infty}$. 
In the following theorems, we present both a point-wise and a uniform in $t$ version of this bound in the asymptotic regime of $T, N \rightarrow \infty$ and under only minimal assumptions on the ground-truth winning probabilities $p_{ij}(t)$'s.

\begin{nthm}\label{thm:orc_ineq_score}
Let $\gamma = \gamma(T,N,p_{\min})$ be the probability that  \cref{cond:nec_suff_bt_1} fails and suppose that the kernel bandwidth is chosen as 
    \begin{align*}
        h = \max\left\{\frac{1}{T^{1+\eta}}, \left(\frac{36 (1-p_{\min}) \log N}{C_s^2 D_m (N-1) T}\right)^{\frac{1}{3}}\right\},
    \end{align*}
     for any $\eta > 0$ and some universal constant $C_s$ depending only to $D_m$, $D_M$, and $W$.
    Then, for each fixed time $t \in (0, 1)$ and sufficiently large $N$ and $T$, 
    \begin{equation} \label{eqn:orc_ineq_score}
    \begin{split}
        & \|\hat{\bbeta}(t) - \bbeta^*(t)\|_\infty \leq 48 M(t) \left( \begin{split}
            & \delta_h(t) + C_s h
        \end{split} \right)
    \end{split}
    \end{equation}
    with probability at least $1-\frac{2}{N} - \gamma$ as long as the right hand side is smaller than $\frac{1}{3}$.
\end{nthm}
Next, we strengthen our previous result, which is valid for each fixed time $t$, to a uniform guarantee over the entire time course.

\begin{nthm}\label{thm:orc_ineq_score_uniform}
Let $\gamma = \gamma(T,N,p_{\min})$ be the probability that  \cref{cond:nec_suff_bt_1} fails and suppose that the kernel bandwidth is chosen as
    \begin{align*}
        h = \max\left\{\frac{1}{T^{1+\eta}}, \left(\frac{36 (1-p_{\min}) \log (NT^{3+3\eta})}{C_s^2 D_m (N-1) T}\right)^{\frac{1}{3}}\right\}
    \end{align*}
    and  that $W$ is $L_W$-Lipschitz.
    Then, for sufficiently large $N$ and $T$,
    \begin{equation} \label{eqn:orc_ineq_score_uniform}
    \begin{split}
        & \sup_{t \in [0,1]}\|\hat{\bbeta}(t) - \bbeta(t)^*\|_\infty \leq 48 \sup_{t \in [0,1]} M(t) 
        \left( \begin{split}
            & \delta_h(t) + C_s h
        \end{split}\right) 
    \end{split}
    \end{equation}
     with probability at least $1-\frac{2 h^3}{N} - \gamma$ as long as the right hand side is smaller than $\frac{1}{3}$.
\end{nthm}

\begin{nrmk}\label{remark:compare.to.Simons.Yao}
    The rate of point-wise convergence for the estimation error implied by the previous result is
    \begin{equation}\label{eqn:conv_rate_params_kern_smooth}
    \begin{split}
        & \begin{split}
            & M(t) \left(\delta_h(t) + \max\left\{\frac{1}{T^{1+\eta}}, \left(\sqrt{\frac{\log N}{NT}}\right)^{\frac{2}{3}}\right\}\right)
        \end{split},
    \end{split}
    \end{equation}
    while the rate for uniform convergence  is
    \begin{equation}\label{eqn:conv_rate_unif_params_kern_smooth}
    \begin{split}
        & \begin{split}
            & M(t) \left(\delta_h(t) + \max\left\{\frac{1}{T^{\frac{1}{2}+\eta}}, \left(\sqrt{\frac{\log (NT^{1+\eta})}{NT}}\right)^{\frac{2}{3}} \right\}\right)
        \end{split}.
    \end{split}
    \end{equation}
    
    Importantly, as we can see in the previous results, the proposed time varying estimator $\hat{\bbeta}(t)$ is consistent only provided that the design regularity parameter $\delta_h(t)$ goes to zero. Of course, if all teams play each other a constant number of times, then $\delta_h(t) = 0$ automatically. In general, however, the impact of the design on the estimation accuracy needs to be assessed on a case-by-case basis.

\end{nrmk}


The rate \eqref{eqn:conv_rate_params_kern_smooth} should be compared with the convergence rate to the true parameters under the static Bradley Terry model, which  \citep{simons1999asymptotics} show to be $O_p(\max_{i,j: i \neq j}\exp(\beta_i^* - \beta_j^*) \sqrt{\log N/NT})$. Thus, not surprisingly, in the more challenging dynamic settings with smoothly varying winning probabilities the estimation accuracy decreases. The exponent of $\frac{2}{3}$ in the rate \eqref{eqn:conv_rate_params_kern_smooth} matches the familiar rate for estimating  H\"older continuous function of order $1$.   


%

From \eqref{eqn:conv_rate_params_kern_smooth} and \eqref{eqn:conv_rate_unif_params_kern_smooth} we observe that the desired oracle property on estimated parameters requires rate constraints on $M(t)$. These constraints appear to be strong assumptions without a direct connection to $p_{ij}(t)$'s in our model-agnostic setting. Instead, we circumvent this issue by introducing a more interpretable condition number $K$ (or $p_{\min}$), dependent on $N, T$, given by
\begin{align*}
    K =& \exp\left(\frac{1}{p_{\min}}\right)
\end{align*}
and proving that for each fixed time $t \in (0, 1)$ our desired oracle property follows from a bound on $M(t)$. 

\begin{nthm}\label{thm:orc_ineq_score_p}
    Under the conditions in \cref{thm:orc_ineq_score}
     \begin{equation} \label{eqn:orc_ineq_score_p}
        \begin{split}
            & \|\hat{\bbeta}(t) - \bbeta^*(t)\|_\infty \leq 72 K \left( \begin{split}
                & \delta_h(t) + C_s h
            \end{split} \right)
        \end{split}
        \end{equation}
    and under the conditions \cref{thm:orc_ineq_score_uniform}
        \begin{equation}\label{eqn:orc_ineq_score_uniform_p}
        \begin{split}
            & \sup_{t \in [0,1]}\|\hat{\bbeta}(t) - \bbeta(t)^*\|_\infty \leq 72 K \sup_{t \in [0,1]}
            \left( \begin{split}
                & \delta_h(t) + C_s h
            \end{split}\right) 
        \end{split}
        \end{equation}
    with probability at least $1 - \frac{2}{N} - \gamma$.
\end{nthm}

\begin{nrmk}\label{rmk:bounded_away_0_pmin}
We note that assuming $p_\text{min}$ to be bounded away from $0$ ensures $\sup_{t \in [0,1]} \|\hat{\bbeta}(t) - \bbeta^*(t)\|_\infty \to 0$, with high probability. This assumption means that no team is uniformly dominated by or dominates others (since this implies $1 - p_\text{min}(t)$ is bounded away from 1). This is a reasonable assumption in real-world data such as sports match histories where teams are screened to be competitive with each other. Therefore it is reasonable to only consider matches between teams that do not have vastly different skills,  which is reflected in winning probabilities that are bounded away from $\theset{0, 1}$.
\end{nrmk}



In summary, our proposed method achieves high-probability oracle bounds on the estimation error 
in our general model-agnostic time-varying setting. We provide the proofs for the stated theorems in Appendix Section~\ref{sec:appdx_pf}.

\section{Experiments}\label{sec:Experiments}
We compare our method with some other methods on synthetic data\footnote{Code available at \url{https://github.com/shamindras/bttv-aistats2020}}. We consider both cases where the pairwise comparison data are generated from the Bradley-Terry model and from a different, nonparametric model. 
\subsection{Bradley-Terry Model as the True Model}\label{subsec:exp_synth_BT}
First we conduct simulation experiments in which the Bradley-Terry model is the true model. Given the number of teams $N$ and the number of time points $M$, the synthetic data generation process is as follows:
\begin{enumerate}
    \item For $i \in [N]$, simulate $\bbeta_i \in \reals^M$ as described below;
    \item For $1\leq i<j\leq N$ and $t \in [M]$, set $n_{ij}{(t)}$ and simulate $X{(t)}$ by
    $
    x_{ij}{(t)} \sim \text{Binom}\Big(n_{ij}{(t)},{1}/\{(1+\exp[\beta_j{(t)}-\beta_i{(t)}]\}\Big)$ and $ 
    x_{ji}{(t)}  = n_{ij}{(t)} - x_{ij}{(t)}$.
\end{enumerate}
For each $i \in [N]$, we generate $\bbeta_i\in \mathbb{R}^M$ from a Gaussian process ${\rm GP}({\mu}_i(t),\sigma_i(t,s))$ as follows:
\begin{enumerate}
    \item Set the values of the mean process ${\mu}_i(t)$ for all $t \in [M]$ and get mean vector $\bm{\mu}_i = (\mu_{i}(1),\ldots,\mu_{i}(M))$;
    \item Set the values of the variance process $\sigma_i(t,s)$ at $(s,t) \in [M]^2$, to derive  $\Sigma_i \in \mathbb{R}^{M\times M}$;
    \item Generate a sample $\bbeta_i$ from ${\rm Normal}(\bm{\mu}_i,\Sigma_i)$.
    
\end{enumerate}

\begin{figure}[htb!]
    \centering
 \includegraphics[width = 0.48\textwidth]{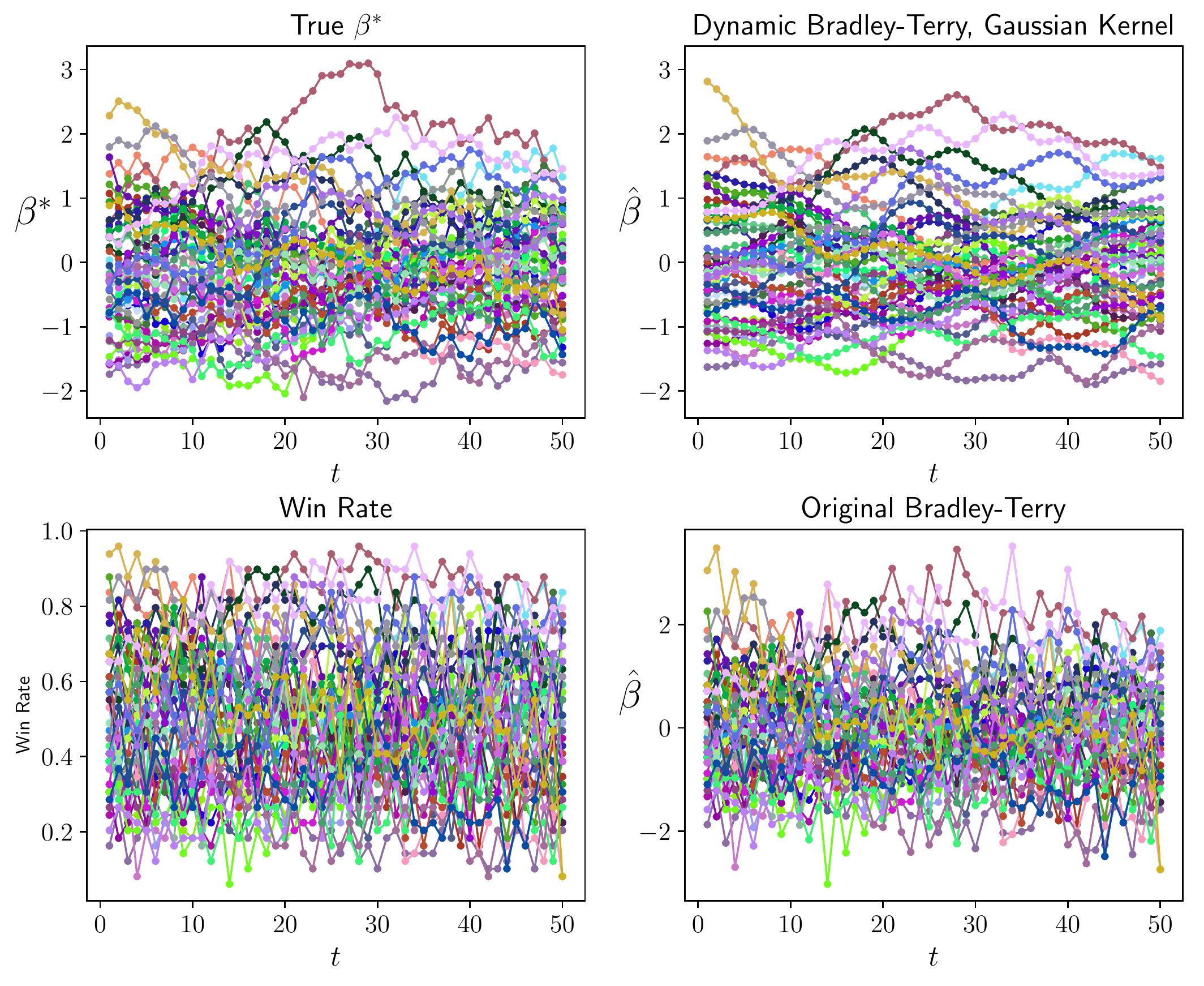}
 \caption{Comparison of $\bbeta$ and different estimators. First row: true $\bbeta$ (left), $\hat{\bbeta}$ with our dynamic BT (right); second row: win rate (left), original BT (right).}
    \label{fig:simulation_result}
\end{figure}

In our experiment, we generate the parameter $\bbeta$ via a Gaussian process for $N = 50$, $M = 50$ and $n_{ij}(t) = 1$ for all $t$. See appendix for full details. We compare the true $\bbeta$, the win rate, and $\hat{\bbeta}$ by different methods in \cref{fig:simulation_result}.

\begin{table}[htb!]
\centering
\scalebox{0.88}{
\begin{tabular}{|c|c|c|c|}
\rowcolor[HTML]{FFFFC7} 
\textbf{Estimator} & \textbf{Rank Diff} & \textbf{LOO Prob} & \textbf{LOO nll} \\
\cellcolor[HTML]{DAE8FC}Win Rate & 3.75 & 0.44 & - \\
\cellcolor[HTML]{DAE8FC}Original BT & 3.75 & 0.37 & 0.56 \\
\cellcolor[HTML]{DAE8FC}Dynamic BT & 2.29 & 0.37 & 0.55
\end{tabular}
}
\caption{Comparison of Different estimators with results based on 20 repeats. Rank Diff. is the average absolute difference between the estimated and true rankings of teams. LOO Prob means average leave-one-out prediction error of win/loss. LOO nll means average leave-one-out negative log-likelihood.}
\label{tab:compare_quant}
\end{table}
 We use LOOCV to select the kernel parameter $h$. The CV curve can be found in \cref{sec:Experiments} of Appendix. In this relatively sparse simulated data we see that our dynamic Bradley-Terry estimator $\hat{\bbeta}$ recovers the comprehensive global ranking of each team in line with the original Bradley-Terry model. We also observe that due to the kernel-smoothing that our estimator has relatively more stable paths over time.
\par
\cref{tab:compare_quant} compares the three estimators across key metrics. The results are averaged over 20 repeats. As expected, in this sparse data setting, our dynamic Bradley-Terry method performs better than the original Bradley-Terry model.

\subsection{Model-Agnostic Setting}\label{subsec:exp_synth_model_free}
\begin{figure}[htb!]
    \centering
 \includegraphics[width = 0.48\textwidth]{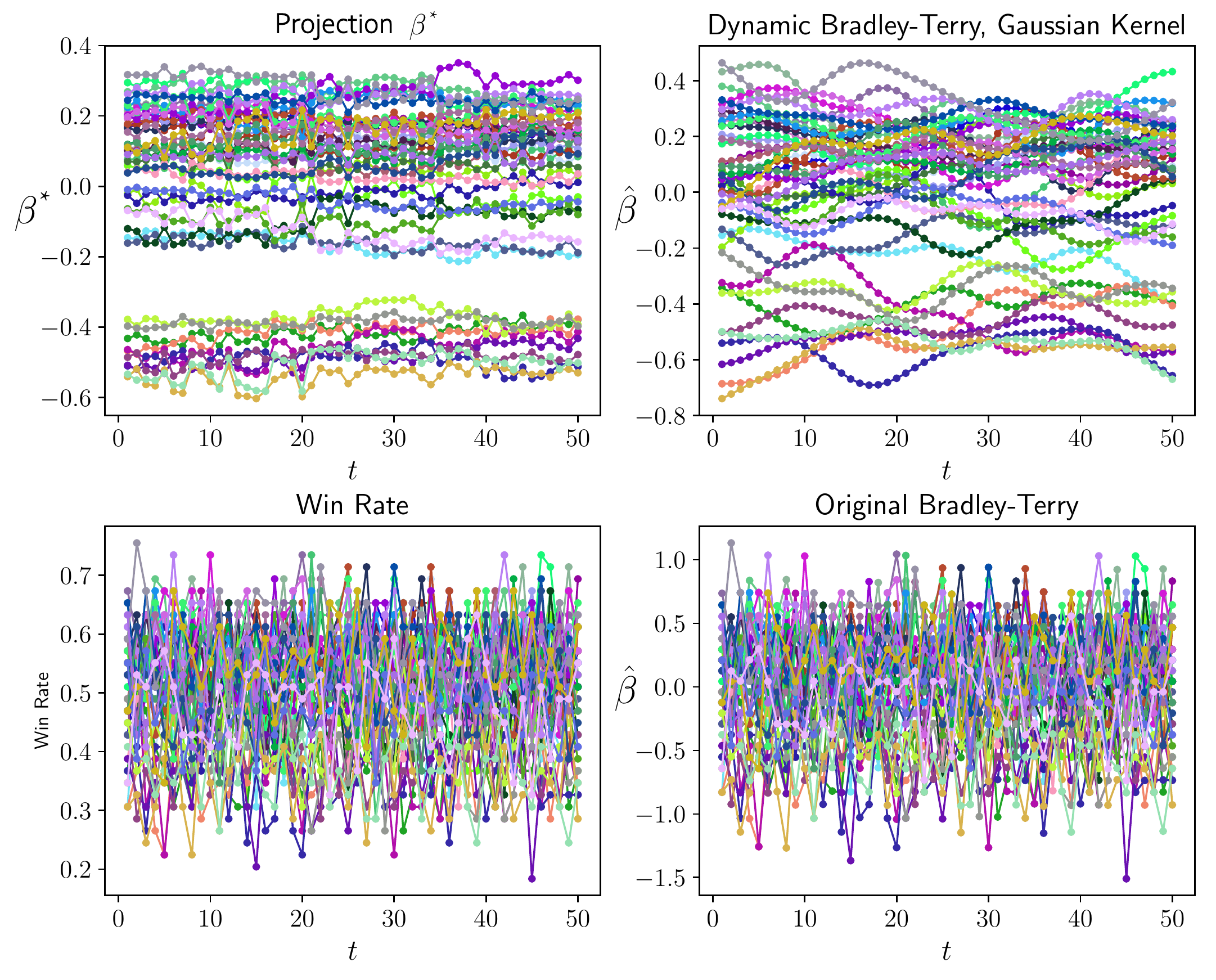}
 \caption{Comparison of $\bbeta^*$ and different estimators when the underlying model is not the Bradley-Terry model. First row: projection $\bbeta^*$ (left), our dynamic BT (right); second row: win rate (left), original BT (right).}
    \label{fig:simulation_result2}
\end{figure}

In our second experiment we adopt a model-agnostic setting where we assume that Bradley-Terry model is not necessarily the true data generating model, as described in \cref{sec:prelim}. With the same notation, for $1\leq i<j\leq N$ and $t\in [M]$ we first simulate $p_{ij}(t)$, and then set $n_{ij}{(t)}$ and simulate $X{(t)}$ by
    $
    x_{ij}{(t)} \sim \text{Binom}\Big(n_{ij}{(t)},p_{ij}(t)\Big)$ and $
    x_{ji}{(t)}  = n_{ij}{(t)} - x_{ij}{(t)}$.
To generate a \textit{smoothly changing} $p_{ij}(t)$, we again use Gaussian process. Specifically, first we generate $p_{ij}(t)$ for $1\leq i<j\leq N$ and $t\in [M]$ from a Gaussian process. Then we scale those $p_{ij}(t)$'s uniformly to make the values fall within a range $[p_l,p_u]$. In our experiment we set $M=50$, $N=50$, $[p_l,p_u] = [0.05,0.95]$ and $n_{ij}(t) = 1$ for all $t$. The projection parameter $\bbeta^*$, the win rate, and $\widehat{\bbeta}$ by different methods are compared in \cref{fig:simulation_result2}. Again the kernel parameter $h$ in our model is selected by LOOCV. 

\begin{table}[htb!]
\centering
\scalebox{0.92}{
\begin{tabular}{|c|c|c|c|}
\rowcolor[HTML]{FFFFC7} 
\textbf{Estimator} & \textbf{Rank Diff} & \textbf{LOO Prob} & \textbf{LOO nll} \\
\cellcolor[HTML]{DAE8FC}Win Rate & 10.68 & 0.49 & - \\
\cellcolor[HTML]{DAE8FC}Original BT & 10.70 & 0.49 & 0.71 \\
\cellcolor[HTML]{DAE8FC}Dynamic BT & 5.48 & 0.49 & 0.68
\end{tabular}
}
\caption{Comparison of different estimators when the underlying model is not the Bradley-Terry model.}
\label{tab:compare_quant2}
\end{table}
 By comparing curves in \cref{fig:simulation_result2}, we note that our estimator $\hat{\bbeta}$ recovers the global rankings better than the win rate and the original Bradley-Terry model, and produces relatively more stable paths over time. The same conclusion is confirmed by \cref{tab:compare_quant2}, which compares the three estimators in some metrics with 20 repetitions.
\par
\begin{nrmk}\label{nrmk:mle-existence}
In this sparse data setting where $n_{ij}(t)$ is fairly small, the original Bradley-Terry model performs worse than our model for two reasons: 1. \cref{cond:nec_suff_bt_1} can fail to hold at some time points, whence the MLE does not exist; 2. even when the MLE exists, it can fluctuate significantly over time because of the relatively small sample size at each time point. As we show in \cref{sec:numer-converge} in the appendix, when $M=50$, $N=50$ and $n_{ij}(t) = 1$, the MLE exists with fairly high frequency. Still, our model performs much better than the original Bradley-Terry model.
\end{nrmk}
\begin{nrmk}\label{nrmk:changepoint-misses}
Since our method is aimed at obtaining accurate estimates of \textit{smoothly} changing beta/rankings with strong prediction power, it may fail to capture some changes in rankings, especially when these changes are relatively small (as in the present case). However the winrate and original Bradley-Terry methods perform much worse, as they appear to miss some true ranking changes while returning many false change points.
\end{nrmk}
\par
Additional details about experiments, including running time efficiency in simulated settings, can be found in the Appendix \cref{sec:exp_synth_ap}.

\begin{table*}[htb!]
    \centering
    \scriptsize
    \begin{tabular}{|c||c|c||c|c||c|c||c|c||c|c|c|c|c|}
    \hline
    \multirow{2}[4]{*}{\textbf{rank}} & \multicolumn{2}{c||}{\textbf{2011}} & \multicolumn{2}{c||}{\textbf{2012}} & \multicolumn{2}{c||}{\textbf{2013}} & \multicolumn{2}{c||}{\textbf{2014}} & \multicolumn{2}{c|}{\textbf{2015}} \bigstrut\\
    \cline{2-11}               & \textbf{ELO}    & \textbf{BT}     & \textbf{ELO}    & \textbf{BT}     & \textbf{ELO}    & \textbf{BT}     &\textbf{ ELO}    & \textbf{BT}     &\textbf{ ELO}    & \textbf{BT} \bigstrut\\
    \hline
        1 &       \cellcolor{blue!25}GB &      \cellcolor{blue!25}GB &       \cellcolor{yellow!25}NE &     \cellcolor{yellow!25}HOU &      \cellcolor{blue!25}SEA &      \cellcolor{blue!25}SEA &      \cellcolor{yellow!25}SEA &     \cellcolor{yellow!25}DEN &      \cellcolor{yellow!25}SEA &     \cellcolor{yellow!25}CAR \\
        2 &       \cellcolor{yellow!25}NE &      \cellcolor{yellow!25}SF &       \cellcolor{yellow!25}DEN &     \cellcolor{yellow!25}ATL &      \cellcolor{yellow!25}SF &      \cellcolor{yellow!25}DEN &      \cellcolor{yellow!25}NE &     \cellcolor{yellow!25}ARI &      \cellcolor{yellow!25}CAR &     \cellcolor{yellow!25}DEN \\
        3 &       \cellcolor{blue!25}NO &      \cellcolor{blue!25}NO &        \cellcolor{yellow!25}GB &     \cellcolor{yellow!25}SF &      \cellcolor{yellow!25}NE &     \cellcolor{yellow!25} NO &      \cellcolor{yellow!25}DEN &     \cellcolor{yellow!25}NE &      \cellcolor{yellow!25}ARI &     \cellcolor{yellow!25}NE \\
        4 &      \cellcolor{yellow!25}PIT &      \cellcolor{yellow!25}NE &        \cellcolor{yellow!25}SF &     CHI &     \cellcolor{yellow!25}DEN &      \cellcolor{yellow!25}KC &      \cellcolor{yellow!25}GB &      \cellcolor{yellow!25}SEA &       \cellcolor{yellow!25}KC &     \cellcolor{yellow!25}CIN \\
        5 &      \cellcolor{yellow!25}BAL &     \cellcolor{yellow!25}DET &       \cellcolor{yellow!25}ATL &     \cellcolor{yellow!25}GB &      \cellcolor{yellow!25}CAR &      \cellcolor{yellow!25}SF &      \cellcolor{blue!25}DAL &     \cellcolor{blue!25}DAL &      \cellcolor{yellow!25}DEN &    \cellcolor{yellow!25}ARI \\
        6 &       \cellcolor{yellow!25}SF &      \cellcolor{yellow!25}BAL &       \cellcolor{yellow!25}SEA &     \cellcolor{yellow!25}NE &      \cellcolor{yellow!25}CIN &      \cellcolor{yellow!25}NE &      PIT &     \cellcolor{yellow!25}GB &      \cellcolor{yellow!25}NE &     \cellcolor{yellow!25}GB \\
        7 &       \cellcolor{yellow!25}ATL &      \cellcolor{yellow!25}PIT &       NYG &     \cellcolor{yellow!25}DEN &     \cellcolor{yellow!25}NO &      \cellcolor{yellow!25}IND &     BAL &      PHI &     \cellcolor{yellow!25}PIT &     \cellcolor{yellow!25}MIN \\
        8 &       PHI &     \cellcolor{yellow!25}HOU &       CIN &     \cellcolor{yellow!25}SEA &     \cellcolor{yellow!25}ARI &   \cellcolor{yellow!25}CAR &   IND &   SD &   \cellcolor{yellow!25}CIN &   \cellcolor{yellow!25}KC\\
        9 &     \cellcolor{yellow!25}SD &    CHI &   \cellcolor{blue!25}BAL &   \cellcolor{blue!25}BAL &   \cellcolor{yellow!25}IND &   \cellcolor{yellow!25}ARI &   \cellcolor{yellow!25}ARI &   DET &   \cellcolor{yellow!25}GB &    \cellcolor{yellow!25}PIT \\
        10 &    \cellcolor{yellow!25}HOU &   \cellcolor{yellow!25}ATL &    \cellcolor{yellow!25}HOU &   IND &   \cellcolor{yellow!25}SD &    \cellcolor{yellow!25}CIN &    CIN &   KC &   \cellcolor{yellow!25}MIN &   \cellcolor{yellow!25}SEA \\
    \hline
            Av. Diff. & \multicolumn{2}{c||}{4.2} & \multicolumn{2}{c||}{5.0} & \multicolumn{2}{c||}{3.5} & \multicolumn{2}{c||}{4.3} & \multicolumn{2}{c|}{3.4} \bigstrut\\
            \hline
    \end{tabular}
    \caption{BT within season vs. ELO NFL top 10 rankings. Blue: perfect match, yellow: top 10 match. Our dynamic BT model is fitted on 16 rounds of each season, and the ranking of a season is based on the ranking at the last round.}
    \label{tab:nfl_elo_topranking}
    
    \end{table*}

\section{Application - NFL Data}\label{sec:application_nfl}
In order to test our model in practical settings we consider
ranking National Football League (NFL) teams over multiple seasons. Specifically we source 5 seasons of openly available NFL data from 2011-2015 (inclusive) using the \texttt{nflWAR} package \citep{yurko2018nflwar}. Each NFL season is comprised of $N = 32$ teams playing $M = 16$ games each over the season. This means that at each point in time $t$ the pairwise comparison matrix based on scores across all 32 teams is sparsely populated with only 16 entries. We fit our time-varying Bradley Terry estimator over all 16 rounds in the season using a standard Gaussian Kernel and tune $h$ using the LOOCV approach described in section \ref{subsec:choosing_lambda_realdata}. In order to gauge whether the rankings produced by our model are reasonable we compare our season-ending rankings (fit over all games played in that season) with the relevant openly available NFL ELO ratings \citep{eloratings2015}. The top 10 season-ending rankings from each method across NFL seasons 2011-2015 are summarized in Table \ref{tab:nfl_elo_topranking}.

Based on Table \ref{tab:nfl_elo_topranking} we observe that we roughly match between 6 to 10 of the top 10 ELO teams consistently over all 5 seasons. There is often misalignment with specific ranking values across both ranking methods. We note that the unlike our method, the NFL ELO rankings use pairwise match data and also additional features including an adjustment for margin of victory. This demonstrates an advantage of our model in only requiring the minimal time-varying pairwise match data and smoothness assumptions to deliver comparable results to this more feature rich ELO ranking method. Furthermore, since our model aggregates data across time it can provide a reasonable minimalist ranking benchmark in modern sparse time-varying data settings with limited ``expert knowledge” e.g. e-sports.

\section{Conclusion}
We propose a time-varying generalization of the Bradley-Terry model that captures temporal dependencies in a nonparametric fashion. This enables the modeling of dynamic global rankings of distinct teams in settings in which the parameteres of the ordinary Bradley Terry model would not be estimable.

From a theoretical perspective we adapt the results of \citep{ford1957solution} to obtain the necessary and sufficient condition for the existence and uniqueness of our Bradley-Terry estimator in the time-varying setting. We extend the previous analysis of \citep{simons1999asymptotics} to derive oracle inequalities on for our proposed method for both the estimation error and the excess risk under smoothness conditions on the winning probabilities. The resulting rates of consistency are of nonparametric type. 

From an implementation perspective we provide a general strategy for tuning the kernel bandwidth hyperparameter using an efficient data-driven approach specific to our unsupervised time-varying setting. Finally, we test the practical effectiveness of our estimator under both simulated and real world settings. In the latter case we separately rank 5 consecutive seasons of open National Football League (NFL) team data \citep{yurko2018nflwar} from 2011-2015. Our NFL ranking results compare favourably to the well-accepted NFL ELO model rankings \citep{eloratings2015}. We thus view our nonparametric time-varying Bradley-Terry estimator as a useful \textit{benchmarking tool} for other feature-rich time-varying ranking models since it simply relies on the minimalist time-varying score information for modeling.





\clearpage
\bibliography{refs}

\clearpage

\section{Appendices}

\subsection{Proofs of Theorems}\label{sec:appdx_pf}

\subsubsection{Proof of Theorem~\ref{thm:exist_unq_thm}}

\paragraph{Uniqueness of the solution}


The elements of the Hessian for $\hat{\mathcal{R}}(\bbeta;t)$ in \eqref{eqn:emp_risk} are:
\begin{equation}
\begin{split}
    H(\hat{\mathcal{R}})_{ii} =& \sum_{j:j \neq i} (\tilde{X}_{ij}(t) + \tilde{X}_{ji}(t)) \frac{\exp{\beta_i}\exp{\beta_j}}{(\exp{\beta_i}+\exp{\beta_j})^2} \\
    H(\hat{\mathcal{R}})_{ij} =& -(\tilde{X}_{ij}(t) + \tilde{X}_{ji}(t)) \frac{\exp{\beta_i}\exp{\beta_j}}{(\exp{\beta_i}+\exp{\beta_j})^2}
\end{split}
\end{equation}

Note that the Hessian has positive diagonal elements, non-positive off-diagonal elements, and zero column sums. With Condition~\ref{cond:nec_suff_bt_1}, this implies that the Hessian can be regarded as a graph Laplacian for a connected graph. Following the classical proof of the property of graph Laplacian \citep{vonLuxburg2007}, 
\begin{equation}
v^TH(\hat{\mathcal{R}})v = \sum_{i<j}\frac{|\tilde{X}_{ij}(t)+\tilde{X}_{ji}(t)|}{2}(v_i - v_j)^2\geq 0
\end{equation}
Then, Condition~\ref{cond:nec_suff_bt_1} guarantees that ``='' is achieved if and only if $v = c\mathbf{1}$. This proves the uniqueness up to constant shifts.


\paragraph{Existence of solution}

Plugging in $\bbeta = \mathbf{0}$, we get an upperbound for the minimum loss function $\hat{\mathcal{R}}^\star(t) := \hat{\mathcal{R}}(\hat{\bbeta};t)$:
\begin{equation}
    \hat{\mathcal{R}}^\star(t) \leq \log 2
\end{equation}
As $\hat{\mathcal{R}}(\bbeta;t)$ is continuous with respect to $\bbeta$, it suffices to show that the level set of $\hat{\mathcal{R}}(\cdot;t)$ at $\log 2$ within $\{\bbeta: \sum_{i=1}^{N}\beta_i = 0\}$ is bounded so that it is compact.

Suppose that $\bbeta$ is in the intersection between the levelset and $\{\bbeta: \sum_{i=1}^{N} \beta_i = 0\}$. Since each summand of $\hat{\mathcal{R}}(\bbeta)$ in \eqref{eqn:emp_risk} is non-negative, i.e.,
\begin{equation}
    \frac{\tilde{X}_{ij}(t)}{\sum_{i',j': i' \neq j'} \tilde{X}_{i'j'}(t)} \log(1 + \exp(\beta_j - \beta_i)) \geq 0
\end{equation} 
if $i$ and $j$ satisfies $\tilde{X}_{ij}(t) > 0$ then the corresponding summand should be smaller than $\log 2$ so that:
\begin{equation}
\begin{split}
    \beta_j - \beta_i \leq& \log(1 + \exp(\beta_j - \beta_i)) \\
    \leq& \log 2 \frac{\sum_{i',j': i'\neq j'} \tilde{X}_{i'j'}(t)}{\tilde{X}_{ij}(t)} 
\end{split}
\end{equation}

By Condition~\ref{cond:nec_suff_bt_1}, for any distinct $i$ and $j$, there exists an index sequence $(i = i_0, i_1, \dots, i_n = j$ such that $X_{i_{k-1} i_k} > 0$ for $k = 1, 2, \dots, n$. Hence,
\begin{equation}
\begin{split}
    \beta_j - \beta_i \leq& \log 2 \sum_{k=1}^{n} \frac{\sum_{i',j': i' \neq j'} \tilde{X}_{i'j'}(t)}{\tilde{X}_{i_{k-1} i_k}(t)}  \\
    \leq& \log2 \sum_{i',j': i' \neq j'} \tilde{X}_{i'j'}(t) \sum_{i',j': i' \neq j'} \frac{1}{\tilde{X}_{i'j'}(t)} =: B
\end{split}
\end{equation}
where $B \in (0, \infty)$.

In sum,
\begin{equation}
    \|\bbeta\|_\infty \leq \max_{i,j: i \neq j} |\beta_i - \beta_j| \leq B
\end{equation}
and this proves the existence part of the theorem.

\subsubsection{Proof of Theorem~\ref{thm:asymp_as_cond}}

The proof of this theorem is based on the proof of Lemma~1 in \cite{simons1999asymptotics}.

Since the kernel function $W$ in Assumption~\ref{assmp:kernel_function} has support $(-\infty, \infty)$, $\tilde{X}_{ij}(t) > 0$ if and only if team $i$ defeated team $j$ at least once any time. In other words, if Condition~\ref{cond:nec_suff_bt_1} holds for at least one time point, then so it does for every time point. Here, we prove that the probability of Condition~\ref{cond:nec_suff_bt_1} to hold at at least one time point converge to $1$ as $N, T \rightarrow \infty$.

Given $p_\text{min}$ instead of $\max_{i,j: i \neq j} \exp(\beta^*_i - \beta^*_j)$, the probability of the event $\mathcal{S}$ that no team in a subset $S$ loses against a team not of $S$ is no larger than
\begin{equation}
     (1-p_\text{min})^{|S| (N - |S| - 1) T}
\end{equation}

Hence, we bound the probability that data does not meet Condition~\ref{cond:nec_suff_bt_1} by a union bound
\begin{equation}
\begin{split}
    & \prob{\text{Condition~\ref{cond:nec_suff_bt_1} fails}} \leq \sum_{S \subset [N]: S \neq \emptyset} \prob{\mathcal{S}} \\
    & \leq \sum_{n=1}^{N-1} \binom{N}{n} (1 - p_\text{min})^{n (N - n - 1) T} \\
    & \leq 2 \sum_{n=1}^{\lceil N/2 \rceil} \binom{N}{n} (1 - p_\text{min})^{n (N - n - 1) T} \\
    & \leq 2 \sum_{n=1}^{\lceil N/2 \rceil} \binom{N}{n} (1 - p_\text{min})^{nNT/2} \\
    & \leq 2\left[(1 + (1 - p_\text{min})^{NT/2})^N - 1\right] \\
    & \leq 2\left[(1 + e^{-\frac{NTp_\text{min}}{2}})^N - 1\right] \\
    & \leq 4 N e^{-\frac{NTp_\text{min}}{2}}
\end{split}
\end{equation}
as long as $e^{-\frac{NTp_\text{min}}{2}} \leq \frac{\log 2}{N}$. We note that $(1+\frac{\log 2}{N})^N \leq e^{\log 2} = 2 \leq 2 \log 2 + 1 = 2 N \frac{\log2}{N} + 1$. Hence,
\begin{equation} \begin{split}
    \prob{\text{Condition~\ref{cond:nec_suff_bt_1} fails}} \leq 4 N e^{-\frac{NTp_\text{min}}{2}}
\end{split} \end{equation}
as long as $N e^{-\frac{NTp_\text{min}}{2}} \leq \log 2$.

Since $N e^{-\frac{NTp_\text{min}}{2}} \geq \log 2$ implies $4 N e^{-\frac{NTp_\text{min}}{2}}$ to be larger than $1$, the probability bound holds for any $N$, $T$, and $p_\text{min}$.


\subsubsection{Proof of Theorem~\ref{thm:orc_ineq_score}}\label{sec:pf_orc_ineq_score}

For readability, in our notation we will omit the dependence on the time point $t$ in the expressions for $\hat{\bbeta}(t)$ and $\bbeta^*(t)$, unless required for clarity.

In our proofs we rely on the results and arguments of \cite{simons1999asymptotics} to demonstrate consistency for the maximum likelihood estimator in the static Bradley-Terry model with an increasing number of parameters. 
In that setting, the authors parametrize the winning probabilities as $p_{i,j} = \frac{u^*_i}{u^*_i + u^*_j}$, where $u^*_i \equiv \exp(\beta^*_i)$, with $\bbeta^* \in \mathbb{R}^N$ such that $\beta^*_1 = 0$. Then, setting $\Delta u_i = \frac{\hat{u}_i - u^*_i}{u^*_i}$, where $\hat{u}_i$ is the MLE of $u^*_i$ (with $\hat{u}_1 = 0$ by convention), it follows from the proof of Lemma 3 of \cite{simons1999asymptotics} that
\begin{equation}
\begin{split}
    & \max_i \frac{\abs{\Delta u_i}}{\abs{\Delta u_i} + 1} \\
    & \leq \frac{8}{N-1} \max_{i,j} \frac{u_i^*}{u_j^*} \max_i \sum_{j: j \neq i} \left\{ \frac{\hat{u}_i}{\hat{u}_i + \hat{u}_j} - \frac{u^*_i}{u^*_i + u^*_j}\right\} \\
\end{split}
\end{equation}
where $u^*_i = \exp( \beta_i^*)$. Next, the authors derived a high probability upper bound on 
\begin{equation}\label{eqn:target_of_prob_bound}
\begin{split}
    & \max_i \sum_{j: j \neq i} \left\{ \frac{\hat{u}_i}{\hat{u}_i + \hat{u}_j} - \frac{u^*_i}{u^*_i + u^*_j}\right\} 
\end{split}
\end{equation}
using the facts that 
\begin{equation}\label{eqn:gradient_real_risk}
    \sum_{j: j \neq i} p_{ij} = \sum_{j: j \neq i} \frac{u^*_i}{u^*_i + u^*_j}
\end{equation}
and 
\begin{equation}\label{eqn:gradient_emp_risk_fixed}
\begin{split}
    \sum_{j: j \neq i} \frac{X_{ij}}{T} = \sum_{j: j \neq i} \frac{\hat{u}_i}{\hat{u}_i + \hat{u}_j},
\end{split}
\end{equation}
 where $X_{ij}$ is the number of matches in which $i$ defeated $j$. The second identity comes from the first order optimality condition of $\hat{\bbeta}$.


In our time-varying setting, however, the subgradient optimality of $\hat{\bbeta}(t)$ for $\hat{\mathcal{R}}(\bbeta;t)$ only imply that, for each $j$,
\begin{equation}\label{eqn:gradient_emp_risk_varying}
    \sum_{j: j \neq i} \tilde{X}_{ij}(t) = \sum_{j: j \neq i} \tilde{T}_{ij}(t) \frac{e^{\hat{\beta}_i}}{e^{\hat{\beta}_j} + e^{\hat{\beta}_i}}.
\end{equation}
Thus, \cref{eqn:gradient_emp_risk_fixed} does not hold in the dynamic setting, due to different $\tilde{X}_{ij}(t) + \tilde{X}_{ji}(t)$ across all  $j \neq i$. Instead, we have that
\begin{equation}
\begin{split}
    & \frac{1}{N-1} \left(\sum_{j: j \neq i} \frac{\tilde{X}_{ij}(t)}{\tilde{T}_{ij}(t)} - \sum_{j: j \neq i} \frac{e^{\hat{\beta}_i}}{e^{\hat{\beta}_j} + e^{\hat{\beta}_i}} \right) \\
    & = \left(\begin{split}
        & \sum_{j: j \neq i} \left( \frac{1}{N-1} - \frac{\tilde{T}_{ij}(t)}{\tilde{T}_i(t)} \right) \frac{\tilde{X}_{ij}(t)}{\tilde{T}_{ij}(t)} \\
        & + \sum_{j: j \neq i} \left(\frac{\tilde{T}_{ij}(t)}{\tilde{T}_i(t)} - \frac{1}{N-1} \right) \frac{e^{\hat{\beta}_i}}{e^{\hat{\beta}_j} + e^{\hat{\beta}_i}}
    \end{split}\right)
\end{split}
\end{equation}

Since $\frac{\tilde{X}_{ij}(t)}{\tilde{T}_{ij}(t)}, \frac{e^{\hat{\beta}_i}}{e^{\hat{\beta}_j} + e^{\hat{\beta}_i}} < 1$, 
\begin{equation}
\begin{split}
    & \left| \begin{split}
        & \frac{1}{N-1} \left(\sum_{j: j \neq i} \frac{\tilde{X}_{ij}(t)}{\tilde{T}_{ij}(t)} - \sum_{j: j \neq i} \frac{e^{\hat{\beta}_i}}{e^{\hat{\beta}_j} + e^{\hat{\beta}_i}}\right)
    \end{split} \right| \\
    & \leq 2 \delta_h(t)
\end{split}
\end{equation}
and
\begin{equation}
\begin{split}
     & \frac{1}{N-1} \sum_{j: j \neq i} \left\{ \frac{e^{\hat{\beta}_i}}{e^{\hat{\beta}_i} + e^{\hat{\beta}_j}} - \frac{e^{\beta^*_i}}{e^{\beta^*_i} + e^{\beta^*_j}}\right\} \\
     & \leq 2 \delta_h(t) + \frac{1}{N-1} \sum_{j: j \neq i} \left\{ \frac{\tilde{X}_{ij}(t)}{\tilde{T}_{ij}(t)} - p_{ij}(t)\right\}.
\end{split}
\end{equation}

To make the bias-variance trade-off due to kernel smoothing more explicit, we decompose the term
\begin{equation}
\begin{split}
     \sum_{j: j \neq i} \left\{ \frac{\tilde{X}_{ij}(t)}{\tilde{T}_{ij}(t)} - p_{ij}(t)\right\}
\end{split}
\end{equation}
as
\begin{equation}\label{eqn:divde_delta_i}
\begin{split}
    & \sum_{j: j \neq i} \left( \begin{split} 
        & \frac{\sum_k W_h(t_k, t) (\mathbf{1}_{ij}(t_k) - p_{ij}(t_k))}{\sum_k W_h(t_k, t)} \\
        \end{split} \right) \\
    & + \sum_{j: j \neq i} \left( \frac{\sum_k W_h(t_k, t) p_{ij}(t_k)}{\sum_k W_h(t_k, t)} - p_{ij}(t) \right) \\
    & =: \Delta^{(var)}_{i} + \Delta^{(bias)}_{i}
\end{split}
\end{equation}
where, for brevity, $t_k$ and $\mathbf{1}_{ij}(t_k)$ here stand for $t_k^{(i,j)}$ and $\mathbf{1}(i \text{ defeats } j \text{ at } t_k)$, respectively.

For the first term, we have that 
\begin{equation}\label{eqn:prob_bound_delta_1}
\begin{split}
    & \prob{\left|\Delta^{(var)}_i\right| \geq \epsilon} \\
    & = \prob{\left|\sum_{j:j \neq i} \frac{\sum_k W_h(t_k, t) (\mathbf{1}_{ij}(t_k) - p_{ij}(t_k))}{\sum_k W_h(t_k, t)}\right| \geq \epsilon} \\
    & = \prob{\begin{split}
        &\left|\sum_{j,k} h W_h(t_k,t) \frac{s_{\min}}{s_j} (\mathbf{1}_{ij}(t_k) - p_{ij}(t_k))\right| \\ 
        & \geq \epsilon \cdot h \cdot s_{\min}
    \end{split}} \\
\end{split}
\end{equation}
where $s_j = \sum_k W_h(t_k,t)$ and $s_{\min} = \min_{j:j\neq i}~s_j$.

Next, $h W_h(t_k,t) \frac{s_{\min}}{s_j} = W\left(\frac{t_k-t}{h}\right) \frac{s_{\min}}{s_j} \leq 1$ and hence that multiplicative Chernoff bound  \citep[see, e.g.][]{raghavan1988probabilistic} yields that
\begin{equation} 
\begin{split}
    & \prob{\left|\Delta^{(var)}_i\right| \geq \epsilon} \\
    & \leq 2\exp\left(- \frac{(\epsilon \cdot h \cdot s_{\min})^2}{3 \sum_{j,k} h W_h(t_k, t) \frac{s_{\min}}{s_j} p_{ij}(t_k)} \right) \\
    & \leq 2\exp\left(- \frac{\epsilon^2 h D_m T}{18 (N-1) (1-p_\text{min})}\right)
\end{split}
\end{equation}
for each $i$ as long as 
\begin{equation} \begin{split}
    \frac{\epsilon}{\sum_{j,k} \frac{W_h(t_k, t)}{\sum_{k'} W_h(t_{k'},t)} p_{ij}(t_k)} \leq 1.
\end{split} \end{equation} 
This condition holds for $\epsilon \leq p_{\min}$.

We note that we have also used the bounds 
\begin{equation}
    \frac{1}{6} D_m T \leq \sum_k W_h(t_k,t) \leq D_M T    
\end{equation} 
for any $i, j$ and sufficiently small $h$, which were shown in \cref{sec:pf_lem_approx_error}.

Then using the union bound,
\begin{equation}
\begin{split}
    & \prob{\max_i \left|\Delta^{(var)}_i\right| \geq \epsilon} \\
    & \leq 2 N \exp\left(- \frac{\epsilon^2 h D_m T}{18 (1-p_\text{min})(N-1) }\right)
\end{split}
\end{equation}
Hence, plugging in $\epsilon = \sqrt{\frac{36 (1 - p_\text{min}) (N-1) \log N}{h D_m T}}$, we get that, with probability at least $1-\frac{2}{N}$, 
\begin{equation}\label{eqn:delta_i_1}
\begin{split}
    \max_i |\Delta^{(var)}_i| \leq& \sqrt{\frac{36 (1 - p_\text{min}) (N-1) \log N}{h D_m T}}
\end{split}
\end{equation}

To handle the deterministic bias terms $\Delta^{(bias)}_{i}$, we rely on the following bound, whose proof is given below in \cref{sec:pf_lem_approx_error}.

\begin{nlem}\label{lem:approx_error}
Suppose that
\begin{enumerate}
    \item $t_1, t_2, \dots, t_T$ satisfies \cref{eqn:bounded_density} and
    \item $\frac{1}{T} = o(h)$ as $T \rightarrow \infty$.
\end{enumerate}
Then, for a $L_f$-Lipschitz function $f: [0,1] \rightarrow \reals$, 
\begin{equation}
\begin{split}
    & \sup_{t \in [0,1]} \left|\sum_{k=1}^{T} \frac{W_h(t_k, t)}{\sum_{k'}W_h(t_{k'},t)} f(t_k) - f(t)\right| \leq C_s h
\end{split}
\end{equation}
with a universal constant $C_s$ depending only on $D_m, D_M, W$ and $L_f$.
\end{nlem}


Accordingly,
\begin{equation}\label{eqn:delta_i_2}
\begin{split}
    & \max_i |\Delta^{(bias)}_i| \\
    & \leq \max_i \sum_{j:j \neq i} \left\{
        \frac{\sum_k W_h(t_k, t)p_{ij}(t_k)}{\sum_k W_h(t_k, t)} - p_{ij}(t)
    \right\} \\
    & \leq C_s (N-1) h
\end{split}
\end{equation}
for some constant $C_s$ depending only on $D_m, D_M, W$ and $L_p$.

Thus, combining all the pieces, 
\begin{equation}
\begin{split}
    & \max_i \frac{|e^{\hat{\beta}_i - \beta^*_i}-1|}{|e^{\hat{\beta}_i - \beta^*_i}-1|+1} \\
    &\leq 8 M(t) \left( 2 \delta_h(t) + \max_i\frac{|\Delta_i^{(var)}(t)| + |\Delta_i^{(bias)}(t)|}{N-1} \right) \\
    & \leq 8 M(t) \left( 2 \delta_h(t) + \sqrt{\frac{36 (1-p_{\min}) \log N}{h D_m (N-1) T}} + C_s h \right)
\end{split}
\end{equation}
with probability at least $1 - \frac{2}{N}$ as long as $\epsilon \leq p_{\min}$.


Plugging in $h = \max\left\{\left(\frac{1}{T}\right)^{1+\eta}, \left(\frac{36 (1-p_{\min}) \log N}{C_s^2 D_m (N-1) T}\right)^{\frac{1}{3}}\right\}$ leads to the bound
\begin{equation}
\begin{split}
    & \max_i \frac{|e^{\hat{\beta}_i - \beta^*_i}-1|}{|e^{\hat{\beta}_i - \beta^*_i}-1|+1} \\
    & \leq 8 M(t) \left( \begin{split}
        & 2 \delta_h(t) + \left(\frac{36 C_s (1-p_{\min}) \log N}{D_m (N-1) T}\right)^{\frac{1}{3}} \\
        & + C_s h
    \end{split} \right) \\
    & \leq 16 M(t) \left( \begin{split}
        & \delta_h(t) + C_s h
    \end{split} \right)
\end{split}
\end{equation}
with probability at least $1-\frac{2}{N}$ when $\epsilon \leq p_{\min}$.
We note that, given our choice for $h$, $\epsilon = \sqrt{\frac{36 (1 - p_\text{min}) (N-1) \log N}{h D_m T}} \leq C_s h$.
Hence, for a sufficiently small $h$, if the right hand side is smaller than, say, $\frac{1}{3}$ then
\begin{equation}
\begin{split}
    & \|\hat{\bbeta} - \bbeta^*\|_\infty \leq 3 \max_i \frac{|e^{\hat{\beta}_i - \beta^*_i}-1|}{|e^{\hat{\beta}_i - \beta^*_i}-1|+1}\\
    & \leq 48 M(t) \left( \begin{split}
        & \delta_h(t) + C_s h
    \end{split} \right)
\end{split}
\end{equation}
with probability at least $1-\frac{2}{N}$ since $\frac{|e^x - 1|}{|e^x - 1|+1} \geq \frac{|x|}{3}$ for $|x| \leq 1$.

\subsubsection{Proof of \cref{lem:approx_error}}\label{sec:pf_lem_approx_error}

Since $f$ is $L_f$-Lipschitz,
\begin{equation}
\begin{split}
    & \left|\sum_{k=1}^{T} \frac{W_h(t_k, t)}{\sum_{k'}W_h(t_{k'},t)} f(t_k) - f(t)\right| \\
    & \leq \sum_{k=1}^{T} \frac{W_h(t_k, t)}{\sum_{k'}W_h(t_{k'},t)} |f(t_k) - f(t)| \\
    & \leq L_f \sum_{k=1}^{T} \frac{W_h(t_k, t)}{\sum_{k'}W_h(t_{k'},t)} |t_k - t|.
\end{split}
\end{equation}

Let $I_1 = \left[0, \frac{t_1+t_2}{2}\right], I_2 = \left[\frac{t_1+t_2}{2}, \frac{t_2+t_3}{2}\right], \dots, I_T = \left[\frac{t_{T-1}+t_T}{2}, 1\right]$ and $l_k$ be the length of $I_k$. We note that $\frac{1}{D_M T} \leq l_k \leq \frac{2}{D_m T}$ by \cref{eqn:bounded_distance}. Then,
\begin{equation}
\begin{split}
    & \int_0^1 |x-t| W_h(x,t) dx = \sum_k \int_{I_k} |x-t| W_h(x,t) \\
    & = \left( \begin{split}
        & \sum_k l_k |t_k-t| W_h(t_k, t) \\
        & + \sum_k \int_{I_k} \left(
        \begin{split}
            & \left|\frac{x-t}{h}\right| W\left(\frac{x-t}{h}\right) \\
            & - \left|\frac{t_k-t}{h}\right| W\left(\frac{t_k-t}{h}\right)
        \end{split}\right)dx \\
    \end{split} \right).
\end{split}
\end{equation}

Since $|\cdot|W$ has a finite total variation,
\begin{equation}
\begin{split}
    & \sum_k \int_{I_k} \left|
        \begin{split}
            & \left|\frac{x-t}{h}\right| W\left(\frac{x-t}{h}\right) \\
            & - \left|\frac{t_k-t}{h}\right| W\left(\frac{t_k-t}{h}\right)
        \end{split}\right|dx \\
    & \leq \sum_k \frac{1}{D_mT} \sup_{x,y \in I_k} \left|
        \begin{split}
            & \left|\frac{x-t}{h}\right| W\left(\frac{x-t}{h}\right) \\
            & - \left|\frac{y-t}{h}\right| W\left(\frac{t_k-t}{h}\right)
        \end{split}\right| \\
    & \leq \frac{\mathcal{V}(|\cdot|W)}{D_mT}.
\end{split}
\end{equation}
Hence,
\begin{equation}
\begin{split}
    & \int_0^1 |x-t| W_h(x,t) dx \\
    & \geq \left( \begin{split}
        & \frac{1}{D_M T} \sum_k |t_k - t| W_h(t_k,t) \\
        & - \frac{\mathcal{V}(|\cdot|W)}{D_mT}
    \end{split} \right).
\end{split}
\end{equation}
As a result,
\begin{equation}\label{eqn:bound_on_sum_xw}
\begin{split}
    & \sum_k |t_k-t| W_h(t_k,t) \\
    & \leq D_M T h \int_{-\infty}^{\infty} |x|W(x)dx + \frac{D_M \mathcal{V}(|\cdot|W)}{D_m}.
\end{split}
\end{equation}

On the other hand, with a similar argument,
\begin{equation}
\begin{split}
    & \int_0^1 W_h(x,t) dx = \sum_k \int_{I_k} W_h(x,t) \\
    & = \left( \begin{split}
        & \sum_k l_k W_h(t_k, t) \\
        & + \sum_k \int_{I_k} \left( \frac{1}{h}W\left(\frac{x-t}{h}\right) - \frac{1}{h}W\left(\frac{t_k-t}{h}\right) \right)dx \\
    \end{split} \right) \\
    & \leq \frac{2}{D_m T} \sum_k W_h(t_k,t) + \frac{\mathcal{V}(W)}{D_m T h},
\end{split}
\end{equation}
implying that
\begin{equation}\label{eqn:bound_on_sum_w}
\begin{split}
    \sum_k W_h(t_k,t) & \geq \frac{D_m T}{2} \int_{-t/h}^{(1-t)/h} W(x)dx - \frac{D_M \mathcal{V}(W)}{2 D_m h}. \\
\end{split}
\end{equation}

As long as $h \rightarrow 0$ and $\frac{1}{T} = o(h)$, 
\begin{equation}
    \inf_{t\in[0,1]} \int_{-t/h}^{(1-t)/h} W(x)dx
\end{equation} 
is bounded below from $0$ (in particular, we consider a small enough $h$ so that it is bounded below by, say, $\frac{1}{3}$), and the term $\frac{D_M \mathcal{V}(|\cdot|W)}{D_m}$ and $\frac{D_M \mathcal{V}(W(x)}{2 D_m h}$ in \cref{eqn:bound_on_sum_xw,eqn:bound_on_sum_w} become asymptotically  negligible. As a result, 
\begin{equation}
\begin{split}
    \sum_{k=1}^{T} \frac{W_h(t_k, t)}{\sum_{k'}W_h(t_{k'},t)} |t_k - t| \leq C' h
\end{split}
\end{equation}
where $C'$ is a universal constant depending only on $D_m$, $D_M$, and $W$, and furthermore 
\begin{equation}
\begin{split}
    \left|\sum_{k=1}^{T} \frac{W_h(t_k, t)}{\sum_{k'}W_h(t_{k'},t)} f(t_k) - f(t)\right| \leq C_s h
\end{split}
\end{equation}
for a univeral constant $C_s$ depending only on $D_m$, $D_M$, $W$, and $L_f$.

\subsubsection{Proof of \cref{thm:orc_ineq_score_uniform}}\label{sec:pf_orc_ineq_score_uniform}

In \cref{sec:pf_orc_ineq_score}, we showed that
\begin{equation}
\begin{split}
    & \max_i \frac{|e^{\hat{\beta}_i(t) - \beta^*_i(t)}-1|}{|e^{\hat{\beta}_i(t) - \beta^*_i(t)}-1|+1} \\
    &\leq 8 M(t) \left( 2 \delta_h(t) + \max_i \frac{|\Delta_i^{(var)}(t)| + |\Delta_i^{(bias)}(t)|}{N-1} \right) \\
\end{split}
\end{equation}

Since the bound for $\Delta^{(bias)}_i(t)$ depends only on $D_m$, $D_M$, $W$, and $L_f$, it is sufficient to find a bound for  
\begin{equation}\label{eqn:sup_max_diff}
    \sup_{t \in [0,1]} \max_i |\Delta^{(var)}_i(t)|.
\end{equation}
We use a covering approach. For $L \in \mathbb{N}$, let $\overline{t}_l = \frac{2l -1}{2L}$ for $l = 1, 2, \dots, L$. Then for any $t \in [0,1]$ there exists $l^*$ such that $|t - \overline{t}_{l^*}| \leq \frac{1}{2L}$ and
\begin{equation}
\begin{split}
    & \max_i \left\{ \Delta^{(var)}_i(t)
    \right\} \\
    & \leq \max_i \left\{ 
    \begin{split} 
    & \sum_{j: j \neq i} \frac{\sum_k W_h(t_k, t) \mathbf{1}(i \text{ defeats } j \text{ at } t_k)}{\sum_k W_h(t_k, t)} \\
    & - \sum_{j: j \neq i} \frac{\sum_k W_h(t_k, \overline{t}_{l^*(t)}) \mathbf{1}(i \text{ defeats } j \text{ at } t_k)}{\sum_k W_h(t_k, \overline{t}_{l^*(t)})}
    \end{split}
    \right\} \\
    & + \max_i \Delta^{(var)}_i(\overline{t}_{l^*})
\end{split}
\end{equation}
where $t_k$ here stands $t_k^{(i,j)}$ for brevity.

In order to bound the second term in the curly brackets, we bound each of its summands as follows:
\begin{equation} \begin{split}
    & \left| \frac{W_h(t_k,t)}{\sum_k W_h(t_k,t)} - \frac{W_h(t_k, \overline{t}_{l^*(t)})}{\sum_k W_h(t_k, \overline{t}_{l^*(t)})} \right| \\
    & \leq \left| \frac{W_t S_{\overline{t}}  - W_{\overline{t}} S_t}{S_t S_{\overline{t}}} \right| \\
    & \leq \frac{|W_t - W_{\overline{t}}|}{S_{\overline{t}}} + \frac{|S_{\overline{t}} - S_t|W_{\overline{t}}}{S_t S_{\overline{t}}} \\
\end{split} \end{equation}
where we denote $W_t = W_h(t_k, t)$, $W_{\overline{t}} = W_h(t_k, \overline{t}_{l^*(t)})$, $S_t = \sum_k W_h(t_k,t)$, and $S_{\overline{t}} = \sum_k W_h(t_k, \overline{t}_{l^*(t)})$ for brevity.

We have seen in \cref{sec:pf_lem_approx_error} that, for any a sufficiently small $h$,
\begin{equation}
    S_t,~S_{\overline{t}} \geq \frac{D_m T}{6}.
\end{equation}

Thus, 
\begin{equation}
\begin{split}
    & \frac{|W_t - W_{\overline{t}}|}{S_{\overline{t}}} + \frac{|S_{\overline{t}} - S_t|W_{\overline{t}}}{S_t S_{\overline{t}}} \\
    & \leq \frac{6}{D_mT} \frac{L_W}{h} |t - \overline{t}_{l^*(t)}| + \left(\frac{6}{D_mT}\right)^2 T \frac{L_W}{h} |t - \overline{t}_{l^*(t)}| \\
    & \leq \frac{36}{D_m^2 L h^2 T}
\end{split}
\end{equation}
as $D_m \leq 1$ and $W$ is $L_W$-Lipschitz by assumption.
Hence, 
\begin{equation}
\begin{split}
    & \max_i \left\{ 
    \begin{split} 
    & \sum_{j: j \neq i} \frac{\sum_k W_h(t_k, t) \mathbf{1}(i \text{ defeats } j \text{ at } t_k)}{\sum_k W_h(t_k, t)} \\
    & - \sum_{j: j \neq i} \frac{\sum_k W_h(t_k, \overline{t}_{l^*(t)}) \mathbf{1}(i \text{ defeats } j \text{ at } t_k)}{\sum_k W_h(t_k, \overline{t}_{l^*(t)})}
    \end{split} \right\} \\
    & \leq \frac{36 (N-1)}{D_m^2 L h^2}
\end{split}
\end{equation}

On the other hand,
\begin{equation}
\begin{split}
    & \max_i \Delta^{(var)}_i(\overline{t}_{l^*}) \leq \max_{l,i} \Delta^{(var)}_i(\overline{t}_l) \\
    & = \max_{l,i} \left\{ 
        \sum_{j: j \neq i} \left( \begin{split} 
        & \frac{\sum_k W_h(t_k, t) (\mathbf{1}_{ij}(t_k) - p_{ij}(t_k))}{\sum_k W_h(t_k, t)} \\
        \end{split} \right)
    \right\}
\end{split}
\end{equation}
where, again,$\mathbf{1}_{ij}(t_k)$ here stands $\mathbf{1}(i \text{ defeats } j \text{ at } t_k)$ for simplicity.

Using \cref{eqn:prob_bound_delta_1} and a union bound, we get that
\begin{equation}
\begin{split}
    & \prob{\max_{l,i} \Delta^{(var)}_i(\overline{t}_l) \geq \epsilon}\\
    & \leq 2NL\exp\left(- \frac{\epsilon^2 h D_m T}{18 (1-p_\text{min}) (N-1)}\right),
\end{split}
\end{equation}
for $\epsilon \leq p_{\min}$.

Next we plug in $\epsilon = \sqrt{\frac{36 (1 - p_{\min}) (N-1) \log(NL)}{h D_m T}}$ to obtain the bounds
\begin{equation}
\begin{split}
    & \left|\max_{l,i} \Delta^{(var)}_i(\overline{t}_{l^*})\right| \leq \sqrt{\frac{36 (1 - p_\text{min}) (N-1) \log (NL)}{h D_m T}}
\end{split}
\end{equation}

and, in turn,
\begin{equation}
\begin{split}
    & \sup_{t,i} \frac{|e^{\hat{\beta}_i(t) - \beta^*_i(t)}-1|}{|e^{\hat{\beta}_i(t) - \beta^*_i(t)}-1|+1} \\
    & \leq \sup_t 8 M(t) \left( \begin{split}
         & 2 \delta_h(t)  + \frac{36}{D_m^2 L h^2}\\
         & + \sqrt{\frac{36 (1-p_{\min}) \log (NL)}{h D_m (N-1) T}} + C_s h
    \end{split} \right)
\end{split}
\end{equation}
with probability at least $1 - \frac{2}{NL}$ and as long as $\epsilon \leq p_{\min}$.

Plugging in $h = \max\left\{\left(\frac{1}{T}\right)^{1+\eta}, \left(\frac{36 (1-p_{\min}) \log (NT^{3+3\eta})}{C_s^2 D_m (N-1) T}\right)^{\frac{1}{3}}\right\}$ and $L = \lceil h^{-3} \rceil$, we conclude that
\begin{equation}
\begin{split}
    & \sup_{t,i} \frac{|e^{\hat{\beta}_i(t) - \beta^*_i(t)}-1|}{|e^{\hat{\beta}_i(t) - \beta^*_i(t)}-1|+1} \\
    & \leq \sup_t 8 M(t) \left( \begin{split}
         & 2 \delta_h(t)  + \frac{36}{D_m^2 L h^2} + C_s h\\
         & + \sqrt{\frac{36 (1-p_{\min}) \log (NL)}{h D_m (N-1) T}}
    \end{split} \right) \\
    & \leq \sup_t 8 M(t) \left( \begin{split}
         & 2 \delta_h(t) + \left(C_s + \frac{72}{D_m^2}\right) h \\
         & + \sqrt{\frac{36 (1-p_{\min}) \log (Nh^{-3})}{h D_m (N-1) T}}
    \end{split} \right) \\
    & \leq \sup_t 16 M(t) \left( \begin{split}
         & \delta_h(t) + \left(C_s + \frac{36}{D_m^2}\right) h\\
    \end{split} \right) \\
\end{split}
\end{equation}
with probability at least $1-\frac{2h^3}{N}$ when $\epsilon \leq p_{\min}$. Since $\epsilon \leq \sqrt{3(1+\eta)} C_s h$ given the choice of  $h$, this bound holds for all sufficiently small $h$.


\subsubsection{Proof of \cref{thm:orc_ineq_score_p}}

For convenience, we omit the time index $t$ for $\hat{\bbeta}(t)$, $\bbeta^*(t)$, and $p_\text{min}(t)$, unless it is required for clarification.

We seek to replace $M(t)$ in \cref{eqn:orc_ineq_score} by a term of $p_\text{min}$. This requires  $\exp(\beta^*_i - \beta^*_j)$ to be bounded above by a function of $p_\text{min}$. The following lemma provides the desired bound. The proof is in \cref{sec:pf_lem_bound_on_beta}.
\begin{nlem}\label{lem:bound_on_beta}
    \begin{equation}
        \max_{i, j: i \neq j} |\beta_i^* - \beta_j^*| - \frac{1}{p_\text{min}}
    \end{equation}
    is upper-bounded by a universal constant, and hence
    \begin{equation}
        \max_{i, j: i \neq j} \exp(|\beta_i^* - \beta_j^*|) \leq C_p \exp \left( \frac{1}{p_\text{min}} \right)
    \end{equation}    
    for some universal constant $1 < C_p < 1.5$.
\end{nlem}

Plugging in the new bound on $\exp(\beta^*_i - \beta^*_j)$, we get 
\begin{equation}
\begin{split}
    & \|\hat{\bbeta}(t) - \bbeta^*(t)\|_\infty \leq 72 K \left( \begin{split}
        & \delta_h(t) + C_s h
    \end{split} \right)
\end{split}
\end{equation}
instead of $48 M(t) \left( \delta_h(t) + C_s h \right)$ in \cref{eqn:orc_ineq_score}

This result easily extends to the uniform case \cref{eqn:orc_ineq_score_uniform_p}.

\subsubsection{Proof of \cref{lem:bound_on_beta}}\label{sec:pf_lem_bound_on_beta}

Let $d_0$ be the difference in scores which implies a bias of probability $\frac{p_\text{min}}{2}$:
\begin{equation}
    \frac{1}{1 + \exp(d_0)} = \frac{p_\text{min}}{2}
\end{equation}
Suppose that 
\begin{equation}
    i_\text{max} = {\arg\max}_i \beta_i^* \text{ and } i_\text{min} = {\arg\min}_i \beta_i^*
\end{equation}
and that
\begin{equation}
    \beta^*_\text{max} = \max_i \beta_i^* \text{ and } \beta^*_\text{min} = \min_i \beta_i^* 
\end{equation}
Then, the maximal difference between $\bbeta_i^*$'s $d_\text{max}$ is 
\begin{equation}
    d_\text{max} = \max_{i,j: i \neq j} \beta_i^* - \beta_j^* = \beta^*_\text{max} - \beta^*_\text{min}
\end{equation}

Let $I_1 = \{i: \beta_i < \beta_\text{min} + d_0\}$. Plugged in $i=i_\text{min}$,  \cref{eqn:gradient_real_risk} implies 
\begin{equation}
\begin{split}
    & (N-1) p_\text{min} \leq \sum_{j:j \neq i_\text{min}} p_{i_\text{min}j}(t) \\
    & = \sum_{j: j \neq i_\text{min}} \frac{1}{1+ \exp(\beta^*_\text{j} - \beta^*_\text{min})} \\
    & \leq \frac{|I_1|-1}{2} + (N-1) \frac{p_\text{min}}{2}
\end{split}
\end{equation}
Hence, $|I_1| \geq (N-1)p_\text{min}+1$. 

Now, let $I_2 = \{i: \beta_i < \beta_\text{min}+2d_0\}$. Summing \cref{eqn:gradient_real_risk} plugged in $i \in I$, we get
\begin{equation}
\begin{split}
    & (N-|I_1|)|I_1|p_\text{min} \leq \sum_{j \in I_1^C} \sum_{i \in I_1} p_{ij}(t) \\
    &= \sum_{j \in I_1^C} \sum_{i \in I_1} \frac{1}{1+ \exp(\beta^*_j - \beta^*_i)} \\
    & \leq \frac{(|I_2|-|I_1|)|I_1|}{2} + (N-|I_1|) |I_1| \frac{p_\text{min}}{2}
\end{split}
\end{equation}
and hence 
\begin{equation}
\begin{split}
    |I_2| \geq& |I_1| + (N-|I_1|)p_\text{min} \\
    \geq& N p_\text{min} + |I_1|(1-p_\text{min}) \\
    \geq& (N-1)(1-(1-p_\text{min})^2) + 1\\
\end{split}
\end{equation}

Similarly for $I_k = \{i: \beta_i < \beta_\text{min} + kd_0\}$ and $J_k = \{j: \beta_j > \beta_\text{max} - kd_0\}$, 
\begin{equation}
\begin{split}
    & |I_k| \geq (N-1)(1-(1-p_\text{min})^k) + 1, \\
    & |J_k| \geq (N-1)(1-(1-p_\text{min})^k) + 1.
\end{split}
\end{equation} 

Now, without loss of generality we assume that $d_\text{max} > 2kd_0$. Then, by the optimality of $\bbeta^*$ for $\mathcal{R}(\bbeta)$,
\begin{equation}
\begin{split}
    \log 2 =& \mathcal{R}(\mathbf{0}) \geq \mathcal{R}(\beta^*) \\
    =& \frac{1}{\binom{N}{2}}\sum_{i,j: i \neq j} p_{ij}(t) \log(1+\exp(\beta_j^* - \beta_i^*)) \\
    \geq& \frac{1}{\binom{N}{2}}\sum_{i \in I_k} \sum_{j \in J_k} p_\text{min} \log(1 + \exp(d_\text{max} - 2k d_0)) \\
    \geq& 2 p_\text{min} (1-(1-p_\text{min})^k)^2 (d_\text{max} - 2kd_0).
\end{split}
\end{equation}
Thus, $d_\text{max} \leq \frac{\log 2}{2 p_\text{min} (1-(1-p_\text{min})^k)^2} + 2kd_0$ for any $k$.
Plugging in $k = \lceil \log(\frac{1}{p_\text{min}})\rceil$, we get that
\begin{equation}
\begin{split}
    d_\text{max} \leq& \frac{\log 2}{2 (1-1/e)^2 p_\text{min}} \\
    & + 2\log\left(\frac{2}{p_\text{min}}-1\right)\left(\log\left(\frac{1}{p_\text{min}}\right)+1\right) \\
    \leq & \frac{1}{p_\text{min}} + C, 
\end{split}
\end{equation}
for some universal constant $C$ since the derivative of $2 \log\left(2x -1 \right)\left(\log x + 1\right)$ is positive and converges to $0$ as $x \rightarrow \infty$. Then $2 \log \left(2x-1\right)\left(\log x + 1\right)$ has a upper-bounding tangent line with slope $1 - \frac{\log 2}{2(1-1/e)^2}$, and $C$ is its y-intercept. This also yields that 
\begin{equation}
    \max_{i,j: i \neq j} \exp(\beta^*_i - \beta^*_j) \leq C_p \exp\left(\frac{1}{p_\text{min}}\right),
\end{equation}
for some universal constant $C_p$. In particular, $1 < C_p < 1.5$.

\subsection{Tuning kernel bandwidth in practical settings}\label{subsec:choosing_lambda_realdata}
As noted in Section \ref{sec:estimator}, the bandwidth $h \in \reals_{> 0}$ of the kernel function serves as an effective global smoothing parameter between subsequent time periods and allows to borrow information across contiguous time points. Increasing $h$, all else held constant, leads to parameter estimates (and hence the derived global rankings) becoming ``smoothed" together across time.

Naturally the question remains on how to \textit{tune} $h$ in practical applications in a principled data-driven manner. This is a fundamentally challenging question not just in our problem but, more generally, in nonparametric regression. Here we present a way to tuning $h$ with some degree of objectivity based on leave-one-out cross-validation (LOOCV).

In general settings where we have independent and identically distributed (i.i.d.) samples, LOOCV assesses the performance of a predictive model on a single held-out i.i.d. sample. In our case, each pairwise comparison can be considered an i.i.d. sample if we take the compared teams and the time point on which they are compared as covariates. Remember that $(i_m, j_m, t_m)$ denotes $m$-th pairwise comparison where team $i_m$ won against team $j_m$ at time point $t_m$ for $m = 1, \dots, M$. Then, for a given smoothing penalty parameter $h$, LOOCV is adapted to our estimation approach as follows:
\begin{enumerate}
    \item For $m = 1, \dots, M$, given $h>0$:
    \begin{enumerate}
        \item fit our model with kernel bandwidth $h$ on the dataset with the $m$-th comparison held-out;
        \item calculate the negative log-likelihood ($\text{nll}$) of the previous solution to $(i_m, j_m, t_m)$.
    \end{enumerate}
    \item Take the average of the negative log-likelihoods to obtain $\text{nll}_{h}$ as a loss in the predictive performance of time-varying Bradley-Terry estimator for given $h$ on our dataset.
    \item Choose the bandwith $h^{*}$ with the smallest $\text{nll}_{h}$ value.
\end{enumerate}
We apply this data-driven methodology to the experiments and real-life application in \cref{sec:Experiments} and \cref{sec:application_nfl}.

\subsection{Details of Experiments} \label{sec:exp_synth_ap}
Here we explain some details of the setting of the numerical experiments in \cref{sec:Experiments}.
\subsubsection{Bradley-Terry Model as the True Model}\label{subsec:exp_synth_BT_ap}
We set the number of teams $N = 50$ and the number of time points $M = 50$. We set $n_{ij}(t) = 1$ for all $i,j\in [N]$ and $t\in [M]$.
\par
For the Gaussian process to generate a path for $\bbeta^{*}_i$ at $t = 1,\ldots,M$, we use the same covariance matrix $\Sigma_i = \Sigma\in \mathbb{R}^{M\times M}$ for all $i \in [N]$, and $\Sigma$ is set to be a Teoplitz matrix defined by
\[
\Sigma_{ij} = 1 - M^{-\alpha}|i-j|^r,
\]
and in our experiment we set $(\alpha,r) = (1,1)$. The mean vector is set to be a constant over time, i.e., $\mu_i(t) = u_i$ for $t = 1,\ldots,M$, and $u_1,\ldots,u_N$ are $i.i.d.$ generated from uniform distribution on $[0,1]$.
 \begin{figure}[htb!]
    \centering
    \includegraphics[width = 0.3\textwidth]{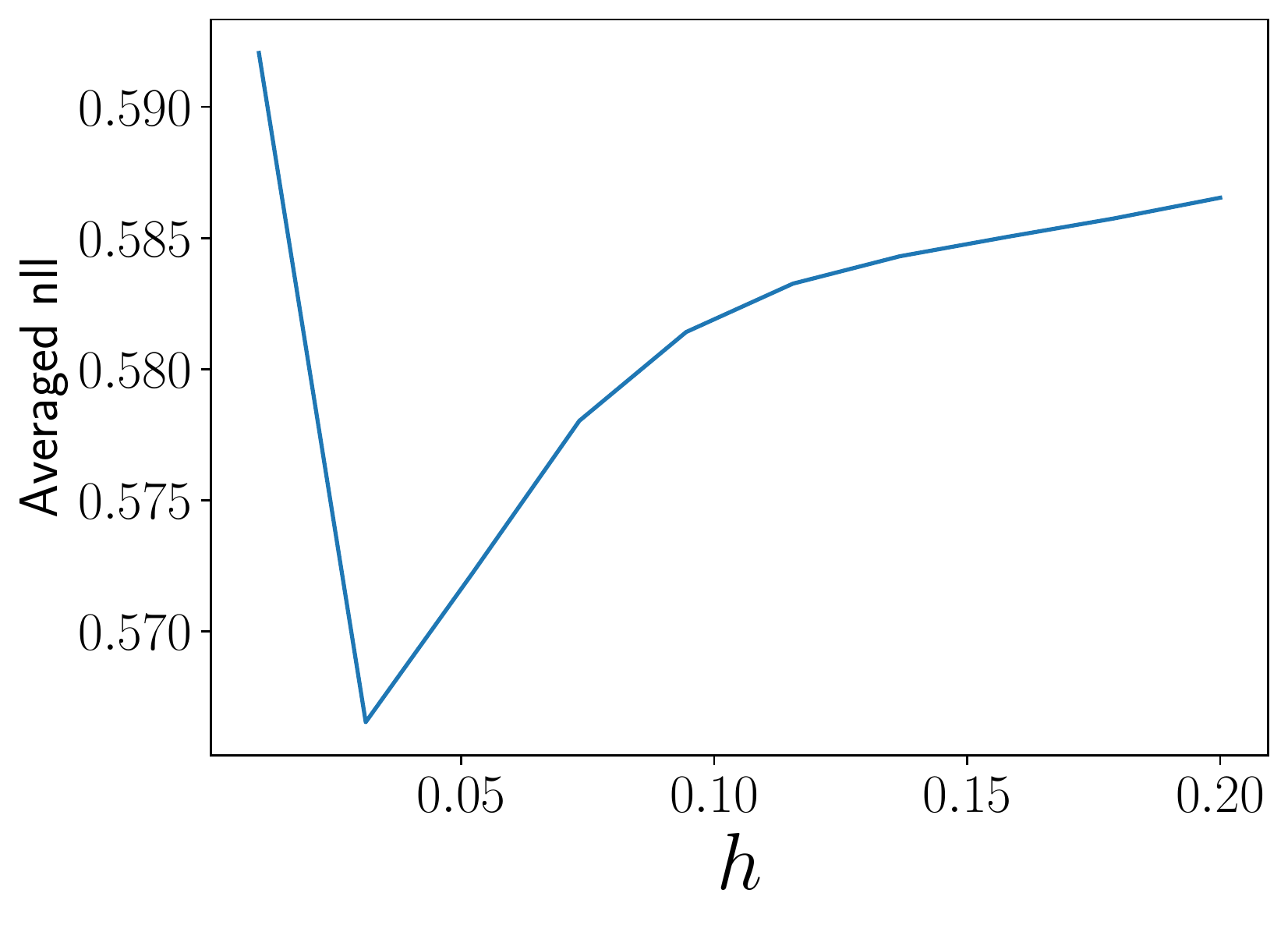}
    \caption{LOOCV curve of our Dynamic Bradley-Terry model fitted with Gaussian kernel. $y$-axis: averaged negative log-likelihood. The optimal $h^*$ is 0.03.}
    \label{fig:cv_curve}
\end{figure}
\par
\cref{fig:cv_curve} shows the curve of LOOCV of our dynamic Bradley-Terry model fitted with a Gaussian kernel in one repetition of our experiment. The curve is for the setting here and for the agnostic model setting the CV curve has similar shape. The curve shows a typical shape of CV curve for tuning parameter. The kernel bandwidth, $h$, with smallest $\mathrm{nll}_{h}$ is $h^{*} = 0.03$. The LOOCV procedure is described in \cref{subsec:choosing_lambda_realdata}.

\subsubsection{Agnostic Model Setting}\label{subsec:exp_synth_model_free_ap}
Again we set the number of teams $N = 50$, the number of time points $M = 50$, and $n_{ij}(t) = 1$ for all $i,j\in [N]$ and $t\in [M]$. The covariance matrix is also the same as section \ref{subsec:exp_synth_BT_ap}. The only difference lies in the mean vector. Now the mean vector is still constant over time, or $\mu_i(t) = u_i$ for $t = 1,\ldots,M$, but $u_i$'s are generated in a following group-wise way:
\begin{enumerate}
    \item Set the number of groups $G$ and the index set of each group $I_1,\ldots,I_G$ so that $\sum_i |I_i| = N$. Set the base support to be $[0,b]$ and the group gap to be $a$.
    \item For each $i\in \{1,\ldots,G\}$, generate $u_j$ from ${\rm Uniform}(a(i-1),a(i-1)+b)$ for all $j\in I_i$. 
\end{enumerate}
In our experiment we set $G = 5$ with each group containing two randomly picked indices, $b = 0.5$ and $a = 1.5$. Such group-wise generation intends to ensure that different teams have distinguishable perofrmance in pairwise comparison so that the ranking is more reasonable.

\subsubsection{Running Time}
\begin{figure}[htb!]
    \centering
 \includegraphics[width = 0.48\textwidth]{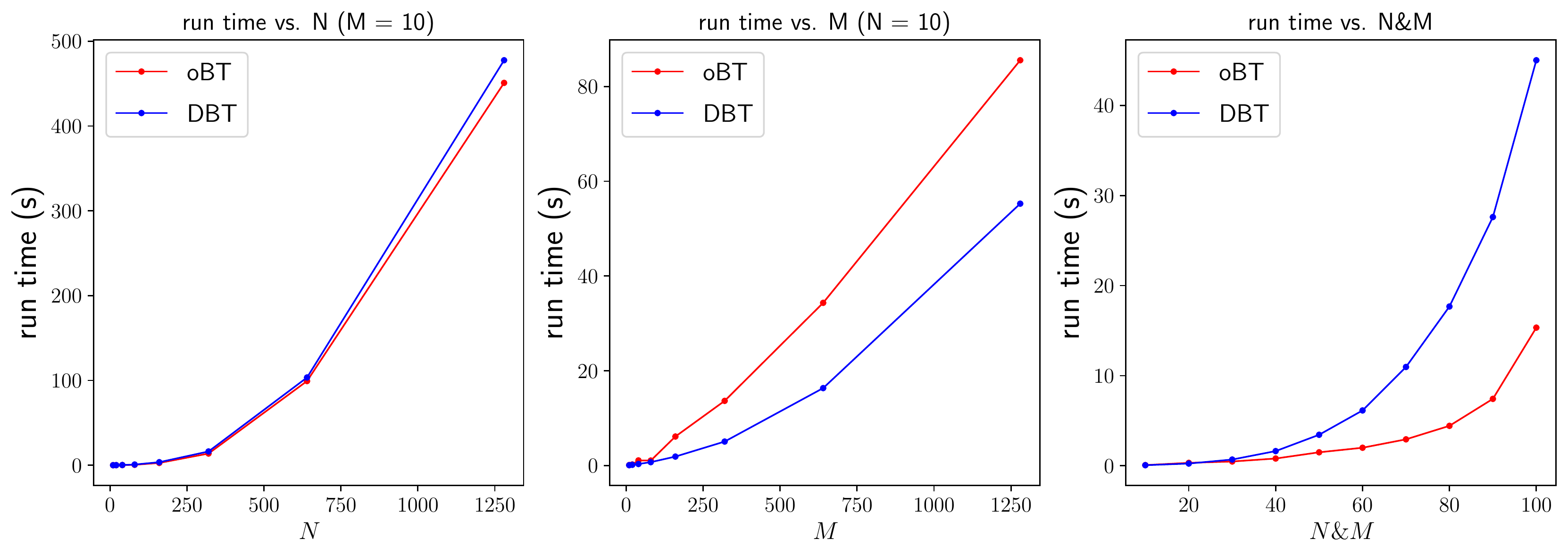}
 \caption{Comparison of running time of original Bradley Terry model (oBT) and our Dynamic Bradley Terry model (DBT). The values are averaged over 20 repetitions.}
    \label{fig:running_time}
\end{figure}

\cref{fig:running_time} compares the time it takes to fit our model and the original Bradley-Terry model under 3 different settings, where $N$ is the number of teams and $M$ is the number of time points:
\begin{itemize}
    \item Fix $N$, vary $M$.
    \item Fix $M$, vary $N$.
    \item Vary $N$ and $M$ together while keeping $N = M$.
\end{itemize}
For our dynamic Bradley-Terry model, the running time here is measured for fitting the model with a given kernel parameter $h$, hence it contains the time cost of kernel smoothing step and the optimization step. In real applications, if one wants to select the best $h$ from a range of values with cross-validation, then the total computation time would be approximately the running time here multiplied by the number of cross-validations.
\par
The results in \cref{fig:running_time} shows that with all advantages our model can bring with, it does not cost much more in terms of computation time. Furthermore, when the number of time points $M$ is large while $N$ is relatively small, our model can cost even less time than the original Bradley-Terry model.
\par
If one wants to do LOOCV to select $h$ when $N$ and $M$ are huge, then it could take a long time to finish the whole procedure. However, in this case we observed in some extended experiments that with a pre-determined $h$ in a reasonable range, our model can give fairly good estimate close to the one given by the best $h^*$ selected by LOOCV. The supporting files of our experiments can be found in our GitHub repository\footnote{Code available at \url{https://github.com/shamindras/bttv-aistats2020}}.
\subsubsection{MLE of the Bradley-Terry Model}
\label{sec:numer-converge}
\cref{tab:connection_original} shows the frequency with which \cref{cond:nec_suff_bt_1} holds at a single time point for the original pairwise comparison data for different $M$ and $N$, where $n_{ij}(t) = 1$ for all $i,j,t$. To be clear, here we just regard the matrix $\tilde{X}(t)$ in \cref{cond:nec_suff_bt_1} as the original data rather than the smoothed data, as it originally was in \cite{ford1957solution}. Given $\{X(t),t\in [M]\}$, the frequency here refers to $\#\{t:\text{The condition holds for } X(t)\}/M$. 
\par
The data are generated as described in \cref{subsec:exp_synth_BT}, and the frequency is averaged over 50 repetitions.  When $N=M=10$ and $n_{ij}(t) = 4$ for all $i,j,t$, the frequency arises to 0.988, illustrating how $n_{ij}(t)$ controls the sparsity of the game matrix and consequently whether \cref{cond:nec_suff_bt_1} holds or not. 
\begin{table}[htb!]
\centering
\scalebox{0.8}{
\begin{tabular}{|c|c|c|c|c|c|c|c|}
$\mathbf{(N,M)}$ & (5,5) & (10,10) & (20,10) & (30,10)& (40,10)& (50,10)\\
Freq. &  0.248 & 0.622 & 0.902 & 0.950& 0.984 & 0.984 \\
\end{tabular}
}
\caption{Frequency that \cref{cond:nec_suff_bt_1} holds at a single time point for the original pairwise comparison data. $n_{ij}(t)=1$.}
\label{tab:connection_original}
\end{table}
\par
As a comparison, under the same setting, for the kernel-smoothed pairwise comparison data, \cref{cond:nec_suff_bt_1} always holds in the experiment. This fact demonstrates the advantage of using kernel-smooth, and partly explains why in our experiments where the data is sparse our model performs the best.
\par
The frequencies in \cref{tab:connection_original} seem high for $N>20$, but from a global perspective, the induced frequency that \cref{cond:nec_suff_bt_1} holds \textit{for all} $M$ time points could be much lower. \cref{tab:connection_original_2} shows such frequency in some settings. Again the values are averaged over 50 repetitions. Remember that in these settings the condition always holds for kernel-smoothed data.
\begin{table}[htb!]
\centering
\scalebox{0.77}{
\begin{tabular}{|c|c|c|c|c|c|c|c|}
$\mathbf{(N,M)}$ & (10,10) &(20,10) & (30,10) & (40,10) & (50,10)& (60,10)\\
Freq. &  0.02 & 0.44 & 0.62 & 0.86& 0.84 & 0.88 \\
\end{tabular}
}
\caption{Frequency that \cref{cond:nec_suff_bt_1} holds for all $M$ time points for the original pairwise comparison data. $n_{ij}(t)=1$.}
\label{tab:connection_original_2}
\end{table}
\par
To make it clearer how $n_{ij}(t)$ affects the global connectivity, we make \cref{tab:connection_original_3}. In the table we fix $(N,M) = (10,10)$.
\begin{table}[htb!]
\centering
\scalebox{0.9}{
\begin{tabular}{|c|c|c|c|c|c|c|c|}
$n_{ij}(t)$ & 1 & 2 & 4 & 6& 8& 10\\
Freq. &  0.02 & 0.48 & 0.92 & 0.94& 0.96 & 1.0 \\
\end{tabular}
}
\caption{Frequency that \cref{cond:nec_suff_bt_1} holds for all $M$ time points for the original pairwise comparison data. $(N,M) = (10,10)$.}
\label{tab:connection_original_3}
\end{table}

\begin{figure}[htb!]
    \centering
 \includegraphics[width = 0.48\textwidth]{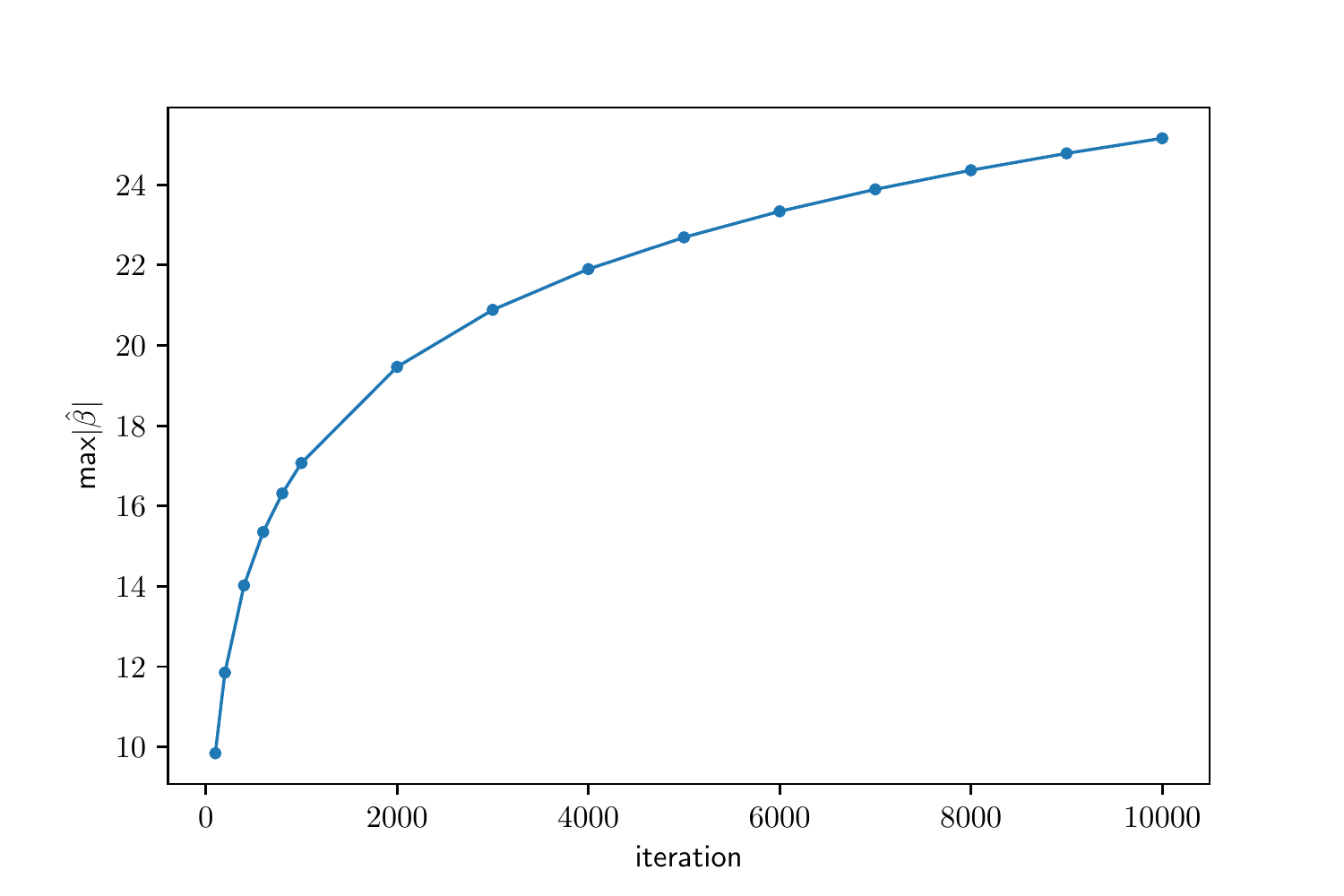}
 \caption{Divergence of max$|\hat{\beta}|$ when MLE does not exist at some time points for the original Bradley-Terry model.}
    \label{fig:diverge}
\end{figure}
\par
By direct inspection of the likelihood  of the original Bradley-Terry model, it can be seen that, when the MLE does not exist, the norm of $\hat{\beta}$ will go to infinity if one uses gradient descent without any regularization. \cref{fig:diverge} shows an example where $N = M =10$ and $n_{ij}(t) = 1$ for all $i,j,t$.

\end{document}